\documentclass[a4paper,11pt,reqno]{amsart}

\usepackage[a4paper,margin=3cm]{geometry}
\usepackage{microtype} 

\usepackage{mathtools,amssymb,amsmath,amsthm}
\usepackage{xfrac}
\usepackage{stmaryrd,mathrsfs,bm,yfonts}
\usepackage{esint}

\usepackage[dvipsnames]{xcolor}
\usepackage[pagebackref,colorlinks=true,linkcolor=MidnightBlue,citecolor=Black,urlcolor=MidnightBlue]{hyperref}
\usepackage[capitalise,nameinlink]{cleveref}
\usepackage{doi}

\usepackage{tikz,tikz-3dplot}

\newtheorem{theorem}{Theorem}[section]
\newtheorem{lemma}[theorem]{Lemma}

\theoremstyle{definition}

\newtheorem{remark}[theorem]{Remark}

\numberwithin{equation}{section}

\renewcommand{\epsilon}{\varepsilon}

\def\XXint#1#2#3{{\setbox0=\hbox{$#1{#2#3}{\int}$ }
  \vcenter{\hbox{$#2#3$ }}\kern-.6\wd0}}

\makeatletter
\renewcommand{\tocsection}[3]{%
  \indentlabel{\@ifnotempty{#2}{\bfseries\ignorespaces\makebox[\@ifempty{#1}{25pt}{75pt}][l]{#1 #2\quad}}}\bfseries#3}
\renewcommand{\tocsubsection}[3]{%
  \indentlabel{\@ifnotempty{#2}{\ignorespaces\makebox[30pt][l]{#1 #2\quad}}}#3}
\newcommand\@dotsep{4}
\def\@tocline#1#2#3#4#5#6#7{\relax
  \ifnum #1>\c@tocdepth \else
    \par \addpenalty\@secpenalty\addvspace{#2}%
    \begingroup \hyphenpenalty\@M
    \@ifempty{#4}{\@tempdima\csname r@tocindent\number#1\endcsname\relax}{\@tempdima#4\relax}%
    \parindent\z@ \leftskip#3\relax \advance\leftskip\@tempdima\relax
    \rightskip\@pnumwidth plus1em \parfillskip-\@pnumwidth
    #5\leavevmode\hskip-\@tempdima{#6}\nobreak
    \leaders\hbox{$\m@th\mkern \@dotsep mu\hbox{.}\mkern \@dotsep mu$}\hfill
    \nobreak
    \hbox to\@pnumwidth{\@tocpagenum{\ifnum#1=1\bfseries\fi#7}}\par
    \nobreak
    \endgroup
  \fi}
\AtBeginDocument{%
\expandafter\renewcommand\csname r@tocindent0\endcsname{0pt}}
\def\l@section{\@tocline{1}{0pt}{1pc}{25pt}{}}
\def\l@subsection{\@tocline{2}{0pt}{2pc}{30pt}{}}
\makeatother

\title[Normalized solution for Kirchhoff equation]{Normalized solution for Kirchhoff equation with upper critical exponent and mixed Choquard type nonlinearities}

\author[J. Shang]{Jinyuan Shang}
\address{School of Mathematics and Computer Sciences, Nanchang University, Nanchang, Jiangxi,\\ 330031, P. R. China}
\email{shangjy0307@163.com}

\author[W. Zhao]{Wenting Zhao}
\address{School of Mathematics and Computer Sciences, Nanchang University, Nanchang, Jiangxi,\\ 330031, P. R. China}
\email{wjkltz@163.com}

\author[X. Huang]{Xianjiu Huang}
\address{School of Mathematics and Computer Sciences, Nanchang University, Nanchang, Jiangxi,\\ 330031, P. R. China}
\email{xjhuangxwen@163.com}
\begin{document}

\maketitle

\begin{abstract}
In this paper, we consider the existence of normalized solution to the following Kirchhoff equation with mixed Choquard type nonlinearities:
\begin{equation*}
\begin{cases}
-\left(a + b \int_{\mathbb{R}^3} |\nabla u|^2 \, dx\right) \Delta u - \lambda u = \mu |u|^{q-2} u + (I_\alpha * |u|^{\alpha + 3}) |u|^{\alpha +1} u, \quad x \in \mathbb{R}^3, \\
\int_{\mathbb{R}^3} u^2 \, dx = \rho^2,
\end{cases}
\end{equation*}
where $a,b,\rho >0$, $\alpha \in \left(0, 3\right)$, $\frac{14}{3} < q < 6$ and $\lambda \in \mathbb{R}$ will arise as a Lagrange multiplier. The quantity $\alpha + 3$ here represents the upper critical exponent relevant to the Hardy-Littlewood-Sobolev inequality, and this exponent can be regarded as equivalent to the Sobolev critical exponent $2^*$. We generalize the results by Wang et al.(Discrete and Continuous Dynamical Systems, 2025), which focused on nonlinear Kirchhoff equations with combined nonlinearities when $2< q< \frac{10}{3}$. The primary challenge lies in the necessity for subtle energy estimates under the \(L^2\)-constraint to achieve compactness recovery. Meanwhile, we need to deal with the difficulties created by the two nonlocal terms appearing in the equation.
\end{abstract}

\renewcommand{\thefootnote}{}
\footnote{Mathematic Subject Classification: 35A15 $\cdot$ 35B33 $\cdot$ 35J60}

\section{Introduction}
In this paper, we study the following nonlinear Kirchhoff equation with mixed Choquard type nonlinearities
\begin{equation}\label{p}
-\left(a + b \int_{\mathbb{R}^3} |\nabla u|^2 \, dx\right) \Delta u - \lambda u = \mu |u|^{q-2} u + (I_\alpha * |u|^{\alpha + 3}) |u|^{\alpha +1} u, \quad x \in \mathbb{R}^3,
\end{equation}
satisfying the mass constraint
\begin{equation}\label{q}
\int_{\mathbb{R}^3} |u|^2 \, dx = \rho^2, \quad u \in H^1(\mathbb{R}^3),
\end{equation}
where $a,b,\rho >0$, $\alpha \in \left(0, 3\right)$, $\frac{14}{3} < q < 6$ and $\lambda \in \mathbb{R}$ will arise as a Lagrange multiplier. $I_\alpha : \mathbb{R}^N \backslash \{0\} \rightarrow \mathbb{R}$ is the Riesz potential of order $\alpha$ defined by
\begin{equation}\label{13}
I_\alpha(x) := \frac{A_\alpha}{|x|^{N-\alpha}} \quad \text{with} \quad A_\alpha = \frac{\Gamma\left(\frac{N-\alpha}{2}\right)}{2^\alpha \pi^{\frac{N}{2}} \Gamma\left(\frac{\alpha}{2}\right)}.
\end{equation}

In the past decade, the nonlinear Kirchhoff equation with combined nonlinearities has attracted extensive attention. The origins of this problem can be traced back to the stationary analogue of the equation
\begin{equation}\label{14}
u_{tt}- \left(a+b\int_{\mathbb{R}^3}|\nabla u|^{2}{d}x \right)\Delta u =f(x,t),
\end{equation}
which was first introduced by Kirchhoff as a generalization of the classical D'Alembert wave equations for the purpose of modeling free vibrations of elastic strings. From a mathematical perspective, (\ref{14}) is widely recognized as a nonlocal equation due to the presence of the term $(\int_{\mathbb{R}^3}|\nabla u|^{2}{d}x)\Delta u$ which prevents (\ref{14}) from being a pointwise identity. Such nonlocality creates major mathematical challenges, rendering the study of Kirchhoff equations both intricate and compelling. Once Lions proposed a functional analysis approach to Kirchhoff equation (\ref{14}) in the pioneer work \cite{Lions}, it captivated significant attention within the mathematical community.

In the equation (\ref{p}), the frequency parameter $\lambda$ can be approached in two distinct ways, giving rise to two separate research directions. When $\lambda$ is prescribed, it is known as the fixed frequency problem, which has been thoroughly researched over the past decade with a wealth of literatures available. The variational method can be utilized to identify the critical points of the energy functional. Alternatively, it is of great interest to investigate solutions to (\ref{p}) that possess a prescribed \(L^2\)-norm (\ref{q}). In this framework, the frequency $\lambda \in \mathbb{R}$ emerges as an unknown, acting as a Lagrange multiplier concerning the constraint (\ref{q}). Consequently, we call it fixed mass problem and the solution is called \textbf{normalized solution}. Physically, identifying normalized solutions is particularly valuable since their \(L^2\)-norm remains constant over time, and analyzing their variational properties can shed light on their orbital stability or instability.

Upon setting the parameters to $a=1$, $b=0$, $\lambda=1$, $\mu=0$, the problem (\ref{p}) is reduced to the celebrated Choquard-Pekar equation
\begin{equation}\label{15}
-\Delta u + u = (I_{\alpha} * |u|^{p})|u|^{p-2}u, \quad x \in \mathbb{R}^N.
\end{equation}
It was Moroz and Van Schaftingen \cite{MVS1} who ascertained an optimal range of the exponent $p$:
\begin{equation*}
\frac{N+\alpha}{N} < p < \frac{N+\alpha}{N-2},
\end{equation*}
within which the existence of solutions to (\ref{15}) is guaranteed when $N \geq 3$. The terms $\frac{N+\alpha}{N}$ and $\frac{N+\alpha}{N-2}$ are respectively designated as the lower and upper critical exponents. These critical exponents are reminiscent of the Sobolev critical exponent in Schr\"{o}dinger equations. For a deeper understanding of the lower critical exponent case, one may refer to \cite{MVS2}, while the upper critical exponent cases are elaborated in \cite{CZ, GY}. For additional related details, the reader is directed to \cite{MVS3, CTW}.

Starting from the foundamental contribution by Tao et al.\cite{TVZ}, Soave \cite{S1} first studied exitence and properties of normalized ground states for the nonlinear Schr\"{o}dinger equation with combined power nonlinearities
\begin{equation}\label{16}
-\Delta u = \lambda u + \mu |u|^{q-2}u + |u|^{p-2}u, \quad x \in \mathbb{R}^N,
\end{equation}
where $2< q \leq 2+\frac{4}{N} \leq p <2^*$. For a $L^2$-subcritical $2 < q < 2 + \frac{4}{N}$, $L^2$-critical $q = 2 + \frac{4}{N}$ and $L^2$-supercritical $2 + \frac{4}{N} < q < 2^*$ perturbation $\mu |u|^{q-2}u$, he proved several existence and stability results in \cite{S1}. Here the Sobolev critical exponent \(2^* := \frac{2N}{N-2}\) for \(N \geq 3\) and \(2^* = \infty\) for \(N = 1, 2\). Later, Soave \cite{S2} considered the Sobolev critical case, i.e., $p=2^*$. By utilizing the fibration method that hinges on the decomposition of the Pohozaev manifold, he explored the existence and nonexistence of normalized solutions to the mixed-dispersion equation (\ref{16}). It should be noted that Jeanjean et al.\cite{JL1} and Wei et al.\cite{Wei} addressed an unresolved issue posed by Soave \cite{S2} in the case where $N \geq 4$ and $N = 3$. They proposed a new set related to the mountain pass level and proved the existence of a second solution of mountain pass type. To build a Palais-Smale sequence at the mountain pass level, the commonly adopted strategy is to utilize the Ghoussoub minimax principle \cite{G1} on the manifold. While this principle served as a very effective theoretical tool for the ranges $2 <q< 2 + \frac{4}{N}$ and $2+\frac{4}{N} \leq q < 2$, it was remarkably intricate because it need to establish much technical topological arguments corresponding to a complex concept named as $\sigma$-homotopy stable family. Thus, a more general minimax principle on the manifold was developed by Chen and Tang \cite{CT1}, which avoids complicated topological arguments, and is technically simpler than \cite{G1} and applicable for all $2 <q< 2^*$. They provided alternative testing functions applicable to all dimensions $N \geq 3$ to estimate the mountain pass level value, thereby achieving the recovery of the compactness of the obtained Palais-Smale sequence. They considered cases where $2 + \frac{4}{N}<q< 2^*$ and $q = 2 + \frac{4}{N}$, as well as the respective dimension $N = 3$ and dimensions $N \geq 4$, within a unified scheme.

Kang and Tang \cite{KT} subsequently genenralized the results of Soave \cite{S2}, they concerned with the nonlinear Schr\"{o}dinger equation with potential and combined power nonlinearities in the Sobolev critical case, where the potential $V$ vanished at infinity.Under certain mild assumptions on $V$, the existence of normalized solutions was established, along with the exponential decay property of positive solutions. This decay property plays a key role in the instability analysis of standing waves. Furthermore, they provided a detailed description of the ground state set and proved the strong instability of standing waves for $q \in [2 + \frac{4}{N}, 2^*)$.

Researchers have focused on the normalized solutions for Choquard equations, considering the following model:
\begin{equation}\label{1a}
-\Delta u = \lambda u + \mu |u|^{q-2}u + (I_\alpha * |u|^p)|u|^{p-2}u, \quad x \in \mathbb{R}^N,
\end{equation}
where $\mu \in \mathbb{R}$ and $\alpha \in (0, N)$, resembling equation (\ref{16}) discussed earlier. Li \cite{L2} demonstrated the existence and orbital stability of ground states of (\ref{1a}) under the conditions where $N\geq 3$, $2 < q < 2+\frac{4}{N}$ and $p=\frac{N+\alpha}{N-2}$. Inspired by Li \cite{L2}, Shang and Ma \cite{SM} investigated the effect of replacing the local perturbation $\mu |u|^{q-2}u$ with the nonlocal term $\mu (I_\alpha *|u|^q)|u|^{q-2}u$ on the existence and behavior of solutions. By varying the parameter $q$, they demonstrated multiple existence theorems based on different conditions imposed on the mass constraint. In addition, Yao et al.\cite{YSW} considered the lower critical exponent case: $p=\frac{N+\alpha}{N}$. Under different assumptions on $q$, $\mu$ and $N \geq 2$, they proved several existence, multiplicity and nonexistence results. Moreover, Giacomoni et al.\cite{GNS} devoted to normalized solution to a critical growth Choquard equation involving mixed operators:
\begin{equation}
-\Delta u + (-\Delta)^s u = \lambda u + g(u) + \left( I_\alpha * |u|^{2^*_\alpha} \right) |u|^{2^*_\alpha - 2}u \quad \text{in} \ \mathbb{R}^N,
\end{equation}
where $2^*_\alpha = \frac{N+\alpha}{N-2}$, $(-\Delta)^s$ is the fractional laplacian operator and $g$ is a real valued function satisfying Ambrosetti-Rabinowitz condition. They looked for a radial solution and investigated the regularity of positive weak solution by using the mountain pass lemma, symmetric decreasing rearrangements and the Strauss lemma.

Recently, Li et al.\cite{LLY} extended the results of Soave, which studied the nonlinear Kirchhoff equations with combined nonlinearities. They considered the existence and asymptotic properties of normalized solutions to the following Kirchhoff equation
\begin{equation}\label{17}
-\left(a + b \int_{\mathbb{R}^3} |\nabla u|^2 \, dx \right) \Delta u = \lambda u + |u|^{p-2}u + \mu |u|^{q-2}u \quad \text{in } \mathbb{R}^3.
\end{equation}
In order to address the unique challenges posed by the nonlocal term $\left(\int_{\mathbb{R}^3} |\nabla u|^2 \, dx\right) \Delta u$  in the Kirchhoff type equations, they developed a perturbed Pohozaev constraint approach and found a way to get a clear picture of the profile of the fiber map via careful analysis. Meanwhile, to recover compactness in the Sobolev critical case, they need some subtle energy estimates under the $L^2$-constraint. Chen and Tang \cite{CT} focused on (\ref{17}) with the Sobolev critical case. The paper \cite{CT} supplies the first non-existence result for $\mu \leq 0$ and $2 < q < 6$. In particular, refined estimates of energy levels are proposed, suggesting a new threshold of compactness in the $L^2$-constraint. This work resolved an open problem within the range $2 < q < \frac{10}{3}$ and completed the analysis for the interval $\frac{10}{3} \leq q < \frac{14}{3}$. It should be noted that the values $\frac{10}{3}$ and $\frac{14}{3}$ correspond to the $L^2$-critical exponents for $N = 3$, as indicated in (\ref{17}) when $b = 0$ and $b > 0$, respectively. From a variational standpoint, besides the Sobolev critical exponent, the \(L^2\)-critical exponent also emerges and plays a central role in the analysis of normalized solutions. This threshold determines whether the constrained functional remains bounded below, thereby influencing the selection of approaches for locating constrained critical points. Their method for constructing Palais-Smale sequences utilized several critical point theorems on manifolds, distinguishing it from the approach in \cite{LLY}. The approach was technically more straightforward and independent of Pohozaev manifold decomposition, making it applicable to a wider class of nonlinear terms exhibiting Sobolev critical growth. Liang \cite{L} considered the above equation (\ref{17}) with well potential and proved the existence of a mountain pass normalized solution via the minimax principle where $\frac{14}{3} < q < p < 6$.

Li et al.\cite{LWM} considered a class nonlocal Kirchhoff-Choquard equation
\begin{equation}\label{18}
-\left(a + b \int_{\mathbb{R}^3} |\nabla u|^2 dx \right) \Delta u - \lambda u = \mu |u|^{q-2}u + \gamma(I_\alpha * |u|^p)|u|^{p-2}u, \quad x \in \mathbb{R}^3
\end{equation}
under different assumptions on $q$, $\mu$, $a$ and $\gamma$, the lower critical exponent case: $p=\frac{N+\alpha}{N}$, they obtained several nonexistence and existence results by new truncation and subadditivity. Wang et al.\cite{WSL} generalized the result of Li et al.\cite{LLY}, which concerned with the qualitative analysis of solutions for the above Kirchhoff equation (\ref{18}) with Choquard mixed nonlinearities, where $2 < q < \frac{10}{3}$, $\frac{14}{3} < p \leq 3 + \alpha$ and $\gamma=1$. When $\frac{14}{3} < p < 3 + \alpha$, they obtained existence and multiplicity of normalized solutions to the above equation, one was a local minimizer and the other of mountain-pass type. When $p=\alpha +3$, they presented a different method of dealing with compactness compared to the general way handling the Sobolev critical term. They applied the Pohozaev identity and $E'_\mu(u_n) \rightarrow 0$ to obtain the sign of the Lagrange multiplier, then the strong convergence of a Palais-Smale sequence in $L^2(\mathbb{R}^3)$ was established. In the defocusing case where $\mu < 0$, an existence result was established through the construction of a minimax characterization for the energy functional. Continuing the exploration initiated by previous researchers, Wang and Sun \cite{WS} further investigated the equation (\ref{18}) under assumptions on \(\mu \in \mathbb{R}\), \(\alpha \in (0, 3)\), \(\frac{10}{3} \leq q< 6\), \(3 + \frac{\alpha}{3} \le p < 3 + \alpha\). They initially considered the case \(\mu > 0\) and identified solutions of mountain pass type. For the defocusing case \(\mu < 0\), the existence of solutions was demonstrated via a minimax characterization of the energy functional. Finally, they discussed the asymptotic behavior of normalized solutions obtained above as \(b \to 0^+\) when \(\mu > 0\).

Inspired by the analysis above, this paper delves into the existence result concerning the Hardy-Littlewood-Sobolev upper critical exponent scenario. It is evident that the energy functional corresponding to problem (\ref{p})-(\ref{q}) is given by
\begin{equation}\label{19}
\Phi(u) := \frac{a}{2}  |\nabla u|^2 _2 + \frac{b}{4}  |\nabla u|^4_2 - \frac{\mu}{q} |u|^q _q - \frac{1}{2(\alpha + 3)} \int_{\mathbb{R}^3} (I_\alpha * |u|^{\alpha + 3}) |u|^{\alpha + 3} \, dx,
\end{equation}
and a critical point of $\Phi$ constrained to
\begin{equation}\label{110}
\mathcal{S}_\rho := \left\{ u \in H^1(\mathbb{R}^3) :  |u|^2_2 = \rho^2 \right\}.
\end{equation}
We recall a solution $u$ to be a ground state solution on $\mathcal{S}_\rho$ if $u$ minimizes the functional $\Phi$ among all the solutions to (1.1), i.e.,
\begin{equation*}
\Phi|'_{\mathcal{S}_\rho}(u) = 0 \text{ and } \Phi(u) = \inf \left\{ \Phi(u) : |u|_2^2 = \rho^2, \ \Phi|'_{\mathcal{S}_\rho}(u) = 0 \right\}.
\end{equation*}

Furthermore, we denote
\begin{equation}\label{111}
S_\alpha := \inf_{u \in D^{1,2}(\mathbb{R}^3) \backslash \{0\}} \frac{\int_{\mathbb{R}^3} |\nabla u|^2 \, dx}{\left( \int_{\mathbb{R}^3} (I_\alpha * |u|^{\alpha +3}) |u|^{\alpha +3} \, dx \right)^{\frac{1}{\alpha +3}}},
\end{equation}
which is similar to the optimal constant of Sobolev inequality. By \cite{SM}, we know that
\begin{equation}\label{112}
S_\alpha = \inf_{u \in D^{1,2}(\mathbb{R}^N) \backslash \{0\}} \frac{|\nabla u|_2^2}{\left( \int_{\mathbb{R}^3} (I_\alpha * |u|^{\alpha + 3}) |u|^{\alpha + 3} \, dx \right)^{\frac{1}{\alpha + 3}}} = \frac{S}{(A_\alpha C_\alpha)^{\frac{1}{\alpha + 3}}},
\end{equation}
where $S$ is the best constant corresponding to the embedding $D^{1,2}(\mathbb{R}^N) \hookrightarrow L^{2^*}(\mathbb{R}^N)$ and $S_\alpha$ is achieved by the family of functions of the form:
\begin{equation*}
U_{\epsilon, x_0}(x) = \frac{(N(N-2)\epsilon^2)^{\frac{N-2}{4}}}{(\epsilon^2 + |x - x_0|^2)^{\frac{N-2}{2}}}, \quad \text{for } x_0 \in \mathbb{R}^N \text{ and } \epsilon > 0,
\end{equation*}
here $C_\alpha = C(N, \alpha)$ given in (\ref{27}) and $A_\alpha$ given in (\ref{13}).

We construct the Pohozaev manifold
\begin{equation}\label{113}
\mathcal{P} := \{ u \in \mathcal{S}_\rho : P(u) = 0 \},
\end{equation}
where
\begin{equation}\label{114}
P(u) := a |\nabla u|^2_2 + b  |\nabla u|^4_2 - \mu \gamma_q  |u|^q _q -\int_{\mathbb{R}^3} (I_\alpha * |u|^{\alpha +3}) |u|^{\alpha +3} \, dx,
\end{equation}
with $\gamma_q = \frac{3(q-2)}{2q}$.

Let $H^1_r(\mathbb{R}^3)$ be the radial subspace of $H^1(\mathbb{R}^3)$. Set
\[
\mathcal{S}_{\rho,r} :=\mathcal{S}_\rho \cap H^1_r(\mathbb{R}^3), \mathcal{P}_r := \mathcal{P} \cap H^1_r(\mathbb{R}^3).
\]
We will show that the energy functional $\Phi$ is coercive on $\mathcal{P}$. Although the embedding of $H^1(\mathbb{R}^3)$ into $L^p(\mathbb{R}^3)$ lacks compactness, the radial subspace $H^1_r(\mathbb{R}^3)$ is compactly embedded in $L^p(\mathbb{R}^3)$ for every $p \in (2, 6)$. Hence, by exploiting symmetric decreasing rearrangements, we seek a radial solution within $H^1_r(\mathbb{R}^3)$. Specifically, the following main result is established::

\begin{theorem}
Assume that $\frac{14}{3} < q < 6$ and $\alpha \in (0,3)$. Then for any $a, b > 0$, $0 <\rho < \kappa$ and $\mu > \mu^*$, where $\kappa$ and $\mu^*$ are given by (\ref{A9}) and (\ref{A10}), there exists a couple $(\lambda_\rho, u_\rho) \in \mathbb{R} \times H^1_r(\mathbb{R}^3)$ solving Eq.(\ref{p})-(\ref{q}) and $\Phi(u) = m_\rho$.
\end{theorem}

\begin{remark}
In comparison with \cite{WSL,WS}, the primary obstacle in the present work lies in the need for delicate energy estimates under the mass constraint to restore compactness. When \( p = 3 + \alpha \), the term \( (I_\alpha * |u|^p)|u|^{p-2}u \) can be seen as a Sobolev critical term, whose presence introduces substantial difficulties due to the lack of compactness. Both Soave \cite{S2} and Li et al. \cite{LLY} applied the approach originally developed by Brezis and Nirenberg \cite{BN} for the Sobolev critical case; this method guarantees that the energy level remains below a specific threshold, which serves as a crucial component in the compactness argument. However, in our article, the existence of \( \left( \int_{\mathbb{R}^3} |\nabla u|^2 dx \right) \Delta u \) and \( (I_\alpha * |u|^p)|u|^{p-2}u \) poses a significant challenge in obtaining precise energy level estimates. Due to the distinctive nature of the upper critical exponent in the Choquard nonlinear term, our method for estimating the upper bound of the mountain pass level differs fundamentally from the approach used in \cite{CT} for the Sobolev critical case. Therefore, establishing the convergence of a Palais-Smale sequence presents a particularly subtle challenge. To overcome these challenges, we follow an approach similar to that in \cite{J} by constructing an auxiliary map on $\mathcal{S}_\rho \times \mathbb{R}$ that retains the same geometric structure as $\Phi(u)$ on $\mathcal{S}_\rho$, and restore compactness via a method analogous to \cite{S1}.
\end{remark}

\begin{remark}
The presence of the term $|\nabla u|_2^4$ in $\Phi(u)$ significantly complicates the establishment of the compactness of Palais-Smale sequences, especially when contrasted with the scenario where $b=0$. This implies that the weak convergence $u_n \rightharpoonup u$ in $H^1(\mathbb{R}^3)$ no longer ensures the convergence
\begin{equation}\label{117}
|\nabla u_n|_2^2 \int_{\mathbb{R}^3} \nabla u_n \cdot \nabla \varphi \, dx \to |\nabla u|_2^2 \int_{\mathbb{R}^3} \nabla u \cdot \nabla \varphi \, dx \quad \text{for all } \varphi \in C_0^\infty(\mathbb{R}^3).
\end{equation}
Consequently, in the case where $b > 0$, excluding the possibilities of vanishing and dichotomy for Palais-Smale sequences becomes increasingly delicate, thereby obstructing their strong convergence in $H^1(\mathbb{R}^3)$. Usually, a bounded Palais-Smale sequence of $\Phi|_{\mathcal{S}_\rho}$ can be obtained by using the Pohozaev constraint approach. That is to say, we can construct a special Palais-Smale sequence $\{u_n\} \subset H^1_{\text{r}}(\mathbb{R}^3)$ for $\Phi|_{\mathcal{S}_\rho}$ with
\[
    P(u_n) = a |\nabla u_n|^2_2 + b  |\nabla u_n|^4_2 - \mu \gamma_q  |u_n|^q _q -\int_{\mathbb{R}^3} (I_\alpha * |u_n|^{\alpha +3}) |u_n|^{\alpha +3} \, dx = o_n(1),
\]
then $\{u_n\}$ is bounded in $H^1(\mathbb{R}^3)$. Once proving $u_n \rightharpoonup u \not\equiv 0$ in $H^1(\mathbb{R}^3)$ for some $u \in H^1(\mathbb{R}^3)$, we can define
\[
    \tilde{B} := \lim_{n \to \infty} |\nabla u_n|_2^2 \geq |\nabla u|_2^2 > 0
\]
and hence (\ref{117}) follows in a standard way. \\
\end{remark}

The paper's structure is as follows: Section 2 delves into an analysis of the initial results, where it is shown that the functional $\Phi$  is coercive on the set $\mathcal{P}_r$. Moreover, it is established that the minimum of $\Phi$ over $\mathcal{P}r$ equals exactly $\sigma\rho$, a strictly positive constant possessing a predetermined upper bound. In Section 3, by analyzing a minimizing sequence of $\Phi$ constrained on $\mathcal{P}_r$, we obtained a weak limit $u\rho$ along with two sequences of Lagrange multipliers ${\lambda_n}$ and ${\zeta_n}$. Utilizing the upper bound for $\sigma_\rho$ and several estimates from Section 2, it is established $u_\rho$ is the strong limit that solves (\ref{p})-(\ref{q}) for $\lambda = \lambda_\rho$.\\

Throughout the paper, we make use of the following notations:
\begin{itemize}
\item $|\cdot|_p$ represents the usual norm of $L^p(\mathbb{R}^3)$. For $1 \leq p < +\infty, L^p(\mathbb{R}^3)$ is the usual Lebesgue space endowed with the norm $|u|_p^p := \int_{\mathbb{R}^3} |u|^p \, dx$. $\|\cdot\|$ represents the usual norm of $H^1(\mathbb{R}^3)$, where $H^1(\mathbb{R}^3)$ is the usual Sobolev space endowed with the norm
    \[
    \|u\|^2 := \int_{\mathbb{R}^3} \left( |\nabla u|^2 + u^2 \right) \, dx.
    \]
    \item $H^1_{r}(\mathbb{R}^3) := \{ u \in H^1(\mathbb{R}^3) \mid u(x) = u(|x|) \text{ a.e. in } \mathbb{R}^3 \}$.
    \item $C_1, C_2, \ldots$ denote positive constants possibly different in different places.
    \item \( o_n(1) \) and \( O_n(1) \) mean that \( |o_n(1)| \to 0 \) and \( |O_n(1)| \leq C \) as \( n \to +\infty \), respectively.
\end{itemize}
\section{Preliminaries}
In this section, several technical results required for the proof of the main theorem are established. We begin by recalling some essential identities.
\begin{lemma}[Gagliardo-Nirenberg inequality, {\cite{W2}}]\label{L6}
Let \( p \in (2, 6) \). Then there exists a constant \( C_p = \left( \frac{p}{2 |\omega_p|_2^{p-2}} \right)^{\frac{1}{p}} > 0 \) such that
\begin{equation}\label{25}
 |u|_p \leq C_p |\nabla u|_2^\delta |u|_2^{1-\delta} \quad \text{for } u \in H^1(\mathbb{R}^3),
\end{equation}
where \( \delta = \frac{3(p-2)}{2p} \) and \( \omega_p \) is the unique positive solution of \( -\Delta \omega + \left( \frac{1}{\delta} - 1 \right) \omega = \frac{2}{p \delta} |\omega|^{p-2} \omega \).
\end{lemma}

\begin{lemma}[Hardy-Littlewood-Sobolev inequality, {\cite{LL}}] \label{L7}
Let $p, r > 1$ and $0 < \alpha < N$ with $\frac{1}{p} + \frac{N-\alpha}{N} + \frac{1}{r} = 2$. Let $u \in L^p(\mathbb{R}^N)$ and $v \in L^r(\mathbb{R}^N)$. Then there exists a sharp constant $C(N, \alpha, p)$ independent of $u$ and $v$ such that
\begin{equation}\label{26}
\left| \int_{\mathbb{R}^N} \int_{\mathbb{R}^N} \frac{u(x) v(y)}{|x-y|^{N-\alpha}} \, dx \, dy \right| \leq C(N, \alpha, p, r) \|u\|_p \|v\|_r.
\end{equation}
If $p = r = \frac{2N}{N+\alpha}$, then
\begin{equation}\label{27}
C(N, \alpha, p, r) := C(N, \alpha) = \pi^{\frac{N-\alpha}{2}} \frac{\Gamma\left(\frac{\alpha}{2}\right)}{\Gamma\left(\frac{N+\alpha}{2}\right)} \left\{ \frac{\Gamma\left(\frac{N}{2}\right)}{\Gamma(N)} \right\}^{-\frac{\alpha}{N}}.
\end{equation}
\end{lemma}
\begin{lemma}\label{L8}{\rm \cite{L1}}
Let $N \geq 3$, $\alpha \in (0, N)$, and $p \in \left[ \frac{N+\alpha}{N}, \frac{N+\alpha}{N-2} \right]$. Assume that the sequence $\{u_n\} \subset H^1(\mathbb{R}^N)$ satisfying $u_n \rightharpoonup u$ in $H^1(\mathbb{R}^N)$ as $n \to \infty$, then
\begin{equation*}
(I_\alpha * |u_n|^p) |u_n|^{p-2} u_n \rightharpoonup (I_\alpha * |u|^p) |u|^{p-2} u \text{ in } H^{-1}(\mathbb{R}^N) \text{ as } n \to \infty.
\end{equation*}
\end{lemma}
\begin{lemma}\label{L9}{\rm \cite{LL}}
Let $u, v, w$ be three Lebesgue measurable non-negative functions on $\mathbb{R}^N$, $N \geq 1$. Then with
\begin{equation*}
\phi(u, v, w) = \int_{\mathbb{R}^N} \int_{\mathbb{R}^N} u(x) v(x-y) w(y) \, dx \, dy,
\end{equation*}
we get
\begin{equation*}
\phi(u, v, w) \leq \phi(u^*, v^*, w^*),
\end{equation*}
where $u^*$ is the Schwartz rearrangement of $u$.
\end{lemma}

Next, we shall introduce the pohozaev manifold and aim to prove that it is a natural constraint.
\begin{lemma}\label{L1}
If $u \in H^1(\mathbb{R}^3)$ is a weak solution of (\ref{p}), then we deduce that $u \in \mathcal{P}$.
\end{lemma}
\begin{proof}
Let $u$ be a solution to (\ref{p}), we have
\begin{equation*}
a |\nabla u|^2_2 + b |\nabla u|^4 _2 = \lambda |u|^2 _2 + \mu |u|^q_q + \int_{\mathbb{R}^3} (I_\alpha * |u|^{\alpha +3}) |u|^{\alpha +3} \, dx.
\end{equation*}
Sceondly, $u$ satisfies the so-called Pohozaev identity
\begin{equation*}
\frac{1}{2} \left( a + b |\nabla u|^2_2 \right) |\nabla u|^2_2 = \frac{3}{2} \lambda |u|^2_2 + \mu \frac{3}{q} |u|^q_q + \frac{1}{2} \int_{\mathbb{R}^3} (I_\alpha * |u|^{\alpha +3}) |u|^{\alpha +3} \, dx.
\end{equation*}
Eliminating the terms with $\lambda$, we have
\begin{equation*}
a |\nabla u|^2_2 + b |\nabla u|^4_2 - \mu \gamma_q |u|^q_q -  \int_{\mathbb{R}^3} (I_\alpha * |u|^{\alpha +3}) |u|^{\alpha +3} \, dx = 0.
\end{equation*}
\end{proof}

Since $\Phi$ is constrained in $\mathcal{S}_\rho$, it is convenient to consider the scaling
\[
t \star u := t^{\frac{3}{2}} u(t \cdot),
\]
which possesses the invariance of $L^2$-norm for $u \in H^1(\mathbb{R}^3) \setminus \{0\}$ and $t > 0$.
we introduce a function
\begin{equation}\label{22}
\aligned
E_u(t) := \Phi(t \star u) =& \frac{a}{2} t^{2} |\nabla u|^2_2 + \frac{b}{4} t^{4} |\nabla u|^4_2 \\
&- \frac{\mu}{q} t^{q \gamma_q } |u|^q_q - \frac{1}{2(\alpha +3)} t^{2(\alpha +3)} \int_{\mathbb{R}^3} (I_\alpha * |u|^{\alpha +3}) |u|^{\alpha +3} \, dx.
\endaligned
\end{equation}
Moreover,
\begin{equation}\label{23}
P(t \star u)=a t^{2} |\nabla u|^2_2 + b t^{4} |\nabla u|^4_2 - \mu \gamma_q t^{q \gamma_q } |u|^q_q -  t^{2(\alpha +3)} \int_{\mathbb{R}^3} (I_\alpha * |u|^{\alpha +3}) |u|^{\alpha +3} \, dx.
\end{equation}

Obviously, we can simply obtain the following result:
\begin{lemma}\label{L2}
For $u \in \mathcal{S}_\rho$, t is a critical point of $E_u(t)$ if and only if $t \star u \in \mathcal{P}$.
\end{lemma}
\begin{proof}
By (\ref{22}), we have
\begin{equation}\label{24}
E'_u(t) = a t |\nabla u|_2^2 + b t^{3} |\nabla u|_2^4 - \mu \gamma_q t^{q \gamma_q -1} |u|_q^q - t^{2\alpha +5} \int_{\mathbb{R}^3} (I_\alpha * |u|^{\alpha +3}) |u|^{\alpha +3} \, dx.
\end{equation}
From the above two identities (\ref{23}) and (\ref{24}), we can get
\begin{equation}\label{F2}
P(t \star u)= tE'_u(t),
\end{equation}
which implies that for any $u \in \mathcal{S}_\rho$, $t \in \mathbb{R}$ is a critical point for $E_u(t)$ if and only if $t \star u \in \mathcal{P}$.
\end{proof}

Next, we consider the decomposition of $\mathcal{P}$ into the disjoint union
\begin{equation*}
\mathcal{P} = \mathcal{P}^+ \cup \mathcal{P}^0 \cup \mathcal{P}^-,
\end{equation*}
where
\begin{equation*}
\mathcal{P}^{+{(\text{resp. } 0, +)} }= \left\{ u \in \mathcal{P} : E_u''(0) > (\text{resp. } =, <) 0 \right\}.
\end{equation*}
\begin{lemma}\label{L3}
Assume that $\frac{14}{3} < q< 6$. Then $\mathcal{P}=\mathcal{P}^-$ in $H^1(\mathbb{R}^3)$.
\end{lemma}
\begin{proof}
For any $u \in \mathcal{P}$, then $P(u) = 0$, namely
\begin{equation*}
a |\nabla u|_2^2 + b |\nabla u|_2^4  - \mu \gamma_q \int_{\mathbb{R}^3} |u|^q \, dx - \int_{\mathbb{R}^3} (I_\alpha * |u|^{\alpha +3}) |u|^{\alpha +3} \, dx = 0.
\end{equation*}
By (\ref{24}), we obtain
\begin{align*}
E''_u(t) =& a |\nabla u|_2^2 + 3b t^{2} |\nabla u|_2^4 - \mu \gamma_q (q\gamma_q-1)t^{q \gamma_q -2} |u|_q^q \\
 &-(2\alpha +5) t^{2\alpha +4} \int_{\mathbb{R}^3} (I_\alpha * |u|^{\alpha +3}) |u|^{\alpha +3} \, dx.
\end{align*}
Therefore, we get
\begin{align*}
E''_u(1) &= a |\nabla u|_2^2 + 3b |\nabla u|_2^4 - \mu \gamma_q (q\gamma_q-1)|u|_q^q -(2\alpha + 5) \int_{\mathbb{R}^3} (I_\alpha * |u|^{\alpha +3}) |u|^{\alpha +3} \, dx \\
&= a |\nabla u|_2^2 + 3b |\nabla u|_2^4 -(2\alpha + 5) \int_{\mathbb{R}^3} (I_\alpha * |u|^{\alpha +3}) |u|^{\alpha +3} \, dx \\
& \quad -(q\gamma_q-1)\left[a |\nabla u|_2^2 + b |\nabla u|_2^4   - \int_{\mathbb{R}^3} (I_\alpha * |u|^{\alpha +3}) |u|^{\alpha +3} \, dx \right]\\
&=(2-q\gamma_q)a |\nabla u|_2^2 + (4-q\gamma_q)b |\nabla u|_2^4 +(q\gamma_q - 2(\alpha +3))\int_{\mathbb{R}^3} (I_\alpha * |u|^{\alpha +3}) |u|^{\alpha +3} \, dx \\
& <0,
\end{align*}
since $\frac{14}{3} < q< 6$ and $\alpha \in (0, 3)$. The last inequality implies that $u \in \mathcal{P}^-$. Hence, we deduce $ \mathcal{P} = \mathcal{P}^-$.
\end{proof}

\begin{lemma}\label{L4}
For every $u \in \mathcal{S}_\rho$, there exists unique $t_u \in \mathbb{R}^+$ such that $t_u \star u \in \mathcal{P}$ where $t_u$ is the strict maxima for $E_u$, that is
\begin{equation*}
\Phi(t_u \star u) = \max_{t > 0} \Phi(t \star u).
\end{equation*}
\end{lemma}
\begin{proof}
From (\ref{22}) and (\ref{24}), we can see $E_u(t)\rightarrow -\infty$ as $t\rightarrow +\infty$ and $E'_u(t)>0$ as $t\rightarrow 0$. Hence, there exists $t_u \in (0, \infty)$ such that $E'_u(t_u) = 0$. Then by Lemma \ref{L2}, $t_u \star u \in \mathcal{P}$. Assume, for contradiction,  that $t_1 \neq t_u$ (take $t_u < t_1$ without loss of generality) satisfies $t_1 \star u \in \mathcal{P}$. As $\mathcal{P} = \mathcal{P}^-$, $t_u$ and $t_1$ are distinct strict local maxima of $E_u$, so there exists $t_2 \in (t_u, t_1)$ for which
\begin{equation*}
E_u(t_2) = \min_{t \in (t_u, t_1)} E_u(t).
\end{equation*}
This leads to $t_2 \in \mathcal{P} \cap \mathcal{P}^+ = \emptyset$, a contradiction. Consequently, there is a unique $t_u \in \mathbb{R}^+$ satisfying $E_u(t_u) = \max_{t > 0} E_u(t)$.
\end{proof}

\begin{remark}\label{L5}
For any $u \in \mathcal{P}$, it follows from (\ref{F2}) that $E_u'(1) = P(1 \star u) = P(u) = 0$. Therefore, $1$ is a critical point of $E_u$. Moreover,  by arguments analogous to Lemma \ref{L4}, $E_u''(1) < 0$, implying that  $t_0 = 1$ is the unique strict local maxima of $E_u$. Hence, for $u \in \mathcal{P}$, it follows that
\begin{equation*}
E_u(1) = \max_{t > 0} E_u(t), \quad \text{or} \quad \Phi(u) = \max_{t > 0} \Phi(t \star u).
\end{equation*}
\end{remark}

\begin{lemma}\label{L10}
For any $u \in \mathcal{P}$, there exists $\delta > 0$ such that  $|\nabla u|_2 \geq \delta$.
\end{lemma}
\begin{proof}
For any $u \in \mathcal{P}$, by Gagliardo-Nirenberg inequality and (\ref{111}), we obtain that
\begin{align*}
a \int_{\mathbb{R}^3} |\nabla u|^2 \, dx + b \left( \int_{\mathbb{R}^3} |\nabla u|^2 \, dx \right)^2 &= \mu \gamma_q \int_{\mathbb{R}^3} |u|^q \, dx +\int_{\mathbb{R}^3} (I_\alpha * |u|^{\alpha +3}) |u|^{\alpha +3} \, dx\\
& \leq  \mu \gamma_q C_q^q \rho^{q(1-\gamma_q)} |\nabla u|_2^{q \gamma_q} + S_\alpha ^{-(\alpha +3)}|\nabla u|_2^{2(\alpha +3)},
\end{align*}
Therefore we have
\begin{equation*}
b |\nabla u|_2^4 \leq \mu \gamma_q C_q^q \rho^{q(1-\gamma_q)} |\nabla u|_2^{q \gamma_q} + S_\alpha ^{-(\alpha +3)}|\nabla u|_2^{2(\alpha +3)}.
\end{equation*}
Since $\frac{14}{3} < q  < 6$ and $\alpha \in (0, 3)$, then $2(\alpha + 3), q \gamma_q > 4$. It follows that
\begin{equation*}
b \leq \mu \gamma_q C_q^q \rho^{q(1-\gamma_q)} |\nabla u|_2^{q \gamma_q - 4} + S_\alpha ^{-(\alpha +3)}|\nabla u|_2^{2(\alpha +1)}.
\end{equation*}
Hence, for any $\mu > 0$, one obtains that there exists $\delta > 0$ satisfying $|\nabla u|_2 \geq \delta$.
\end{proof}

\begin{lemma}\label{L11}
If $u \in \mathcal{P}$ is a critical point for $\Phi|_{\mathcal{P}}$ with $E''_u(1)\neq 0$, then $u$ is a critical point for $\Phi|_{S_\rho}$.
\end{lemma}
\begin{proof}
If $u \in \mathcal{P}$ is a critical point for $\Phi|_{\mathcal{P}}$, by the Lagrange multipliers rule, there exist $\lambda, \zeta \in \mathbb{R}$ such that
\begin{equation}\label{F1}
\Phi'(u) - \lambda u + \zeta P'(u) = 0 \quad \text{in } \mathbb{R}^3,
\end{equation}
where
\[
\Phi'(u) = -\left(a + b |\nabla u|^2_2 \right) \Delta u - \mu |u|^{q-2} u - (I_\alpha * |u|^{\alpha +3}) |u|^{\alpha +1} u,
\]
and
\[
P'(u) = -\left(2a + 4b |\nabla u|^2_2 \right) \Delta u - \mu q\gamma_q |u|^{q-2} u - 2(\alpha +3)(I_\alpha * |u|^{\alpha +3}) |u|^{\alpha +1} u.
\]
It suffices to show that $\zeta = 0$. Observe that any solution $u$ of (\ref{F1}) must satisfy the corresponding Pohozaev identity
\[
(J_u)'(1)= \frac{d J(t^{\frac{3}{2}} u(tx))}{dt} \bigg|_{t=1}=0,
\]
where $J$ is the energy functional corresponding to (\ref{F1}) given by
\[
J(u) := \Phi(u) - \frac{\lambda}{2} |u|_2^2 + \zeta P(u).
\]
Therefore, by using (\ref{F2}) we get
\[
0 = \frac{d J(t^{\frac{3}{2}} u(tx))}{dt} \bigg|_{t=1} = \frac{d J(t \star u)}{dt} \bigg|_{t=1} = E_u'(1) + \zeta (E_u''(1) + E_u'(1)) = \zeta E_u''(1).
\]
which implies that $\zeta = 0$ since $E''_u(1)\neq 0$, so the results follows.
\end{proof}

\begin{lemma}\label{L12}
For $\frac{14}{3} < q < 6$, there is a constant $k > 0$ small enough such that
\[
u \in \overline{A_k} \Rightarrow \Phi(u), \, P(u) > 0,
\]
where $A_k = \{ u \in \mathcal{S}_\rho : |\nabla u|_2^2 < k \}$.
\end{lemma}
\begin{proof}
It follows from (\ref{19}), (\ref{114}), (\ref{111}) and Gagliardo-Nirenberg inequality that
\begin{equation}
\aligned
\Phi (u) &=\frac{a}{2} |\nabla u|_2^2 + \frac{b}{4} |\nabla u|_2^4  - \frac{\mu}{q}  \int_{\mathbb{R}^3} |u|^q \, dx - \frac{1}{2(\alpha +3)} \int_{\mathbb{R}^3} (I_\alpha * |u|^{\alpha +3}) |u|^{\alpha +3} \, dx \\
&\geq \frac{a}{2} |\nabla u|_2^2 + \frac{b}{4} |\nabla u|_2^4 - \frac{\mu}{q} \gamma_q C_q^q \rho^{q(1-\gamma_q)} |\nabla u|_2^{q \gamma_q } - \frac{S_\alpha ^{-(\alpha +3)}}{2(\alpha +3)}|\nabla u|_2^{2 (\alpha +3)},
\endaligned
\end{equation}
and
\begin{equation}
\aligned
P(u) &= a |\nabla u|_2^2 + b |\nabla u|_2^4  - \mu \gamma_q \int_{\mathbb{R}^3} |u|^q \, dx - \int_{\mathbb{R}^3} (I_\alpha * |u|^{\alpha +3}) |u|^{\alpha +3} \, dx \\
&\geq a |\nabla u|_2^2 + b |\nabla u|_2^4 - \mu \gamma_q C_q^q \rho^{q(1-\gamma_q)} |\nabla u|_2^{q \gamma_q } - S_\alpha ^{-(\alpha +3)}|\nabla u|_2^{2(\alpha +3)},
\endaligned
\end{equation}
then $\Phi(u) , P(u) > 0$ for any $u \in \overline{A_k}$ with $k > 0$ small enough. So we complete the proof.
\end{proof}

Define a set of paths
\[
\Gamma_\rho = \left\{ \eta \in C([0,1], \mathcal{S}_{\rho , r} \times \mathbb{R}) : \eta(0) \in \overline{A_k} \text{ and } \eta(1) \in \Phi^0 \right\},
\]
and a min-max value
\[
\sigma_\rho := \inf_{\eta \in \Gamma_\rho} \left\{ \max_{t \in [0,1]} \Phi(\eta(t(x))) \right\}.
\]

We introduce an auxiliary functional as in \cite{J}:
\[
\widetilde{\Phi}: S_{\rho} \times \mathbb{R} \rightarrow \mathbb{R}, \quad (u, t) \mapsto \Phi(t \star u),
\]
to establish a relation between $\sigma_\rho$, stated in the following Lemma, and the infimum of $\Phi$ on $\mathcal{P}$ and $\mathcal{P}_r$. For the same, let us define $\Phi^k := \{ u \in \mathcal{S}_{\rho , r} : \Phi(u) \leq k \}$.
\begin{lemma}\label{L13} {\rm ({\cite{J}} Proposition 2.1)}
Define a set of paths
\[
\tilde{\Gamma}_\rho = \left\{ \tilde{\eta} \in C([0,1], \mathcal{S}_{\rho , r} \times \mathbb{R}) : \tilde{\eta}(0) \in (\overline{A_k}, 1) \text{ and } \tilde{\eta}(1) \in (\Phi^0, 1) \right\}
\]
with associated min-max value
\[
\tilde{\sigma}_\rho := \inf_{\tilde{\eta} \in \tilde{\Gamma}_\rho} \left\{ \max_{t \in [0,1]} \widetilde{\Phi}(\tilde{\eta}(t)) \right\}.
\]
Then, $\tilde{\sigma}_\rho = \sigma_\rho$.
\end{lemma}

In fact, this identity can be directly derived from the definitions of $\widetilde{\sigma}_{\rho}$ and $\sigma_{\rho}$, along with the properties of the maps
\[
\varphi: \Gamma_{\rho} \rightarrow \widetilde{\Gamma}_{\rho}, \eta \rightarrow \varphi(\eta) := (\eta, 0)
\]
and
\[
\psi: \widetilde{\Gamma}_{\rho} \rightarrow \Gamma_{\rho}, \widetilde{\eta} \rightarrow \psi(\widetilde{\eta}) := H \circ \widetilde{\eta}
\]
satisfy
\[
\widetilde{\Phi}(\varphi(\eta)) = \Phi(\eta) \text{ and } \Phi(\psi(\widetilde{\eta})) = \widetilde{\Phi}(\widetilde{\eta}).
\]
where $H(t, u) = e^{\frac{3}{2} t} u(e^{t} x)$.

\begin{lemma}\label{L14} {\rm ({\cite{J}} Proposition 2.2)}
Let $\{h_n\} \subset \tilde{\Gamma}_\rho$ be such that
\[
\max_{x \in [0,1]} \tilde{\Phi}(h_n(x)) \leq \tilde{\sigma}_\rho + \frac{1}{n}, \quad \text{for every } n \in \mathbb{N}.
\]
Then there exists a sequence $\{(u_n, t_n)\} \subset S_{\rho , r} \times \mathbb{R}$ satisfying
\begin{enumerate}
    \item $\tilde{\sigma}_\rho - \frac{1}{n} \leq \tilde{\Phi}(u_n, t_n)\leq \tilde{\sigma}_\rho + \frac{1}{n}$,
    \item $\min_{x \in [0,1]} \|(u_n, t_n) - h_n(x)\| \leq \frac{1}{\sqrt{n}}$,
    \item $\|\tilde{\Phi}'_{S_{\rho , r} \times \mathbb{R}}(u_n, t_n)\| \leq \frac{2}{\sqrt{n}}$, that is
    \[
    |\tilde{\Phi}'(u_n, t_n)(z)| \leq \frac{2}{\sqrt{n}} \|z\| \quad \forall z \in \tilde{\Gamma}_{(u_n, t_n)} = \left\{ (z_1, z_2) \in H^1(\mathbb{R}^N) : \int_{\mathbb{R}^N} u_n z_1 = 0 \right\}.
    \]
\end{enumerate}
\end{lemma}

\begin{lemma}\label{L15}
There exists a sequence $\{v_n\} \subset \mathcal{S}_{\rho , r}$ such that, as $n\rightarrow \infty$
\[
    \Phi(v_n) \rightarrow \sigma_\rho, \quad \Phi|'_{\mathcal{S}_{\rho , r}} \rightarrow 0,
    \quad P(v_n) \rightarrow 0.
\]
\end{lemma}
\begin{proof}
Let $\{\tilde{f}_n\} = \{(f_n, 1)\} \subset \tilde{\Gamma}_\rho$ be such that
\[
\max_{x \in [0,1]} \tilde{\Phi}(\tilde{f}_n(x)) \leq \tilde{\sigma}_\rho + \frac{1}{n} \quad \text{for every } n \in \mathbb{N}.
\]
Then, according to Lemma \ref{L14}, there exists a sequence $\{(u_n, t_n)\} \subset S_{\rho , r} \times \mathbb{R}$ satisfying
\begin{enumerate}
    \item $\tilde{\sigma}_\rho - \frac{1}{n} \leq \tilde{\Phi}(u_n, t_n) \leq \tilde{\sigma}_\rho + \frac{1}{n}$, and hence
    \begin{equation}\label{D1}
    \{\tilde{\Phi}(u_n, t_n)\} \to \tilde{\sigma}_\rho \quad \text{as } n \to \infty,
    \end{equation}
    \item $\min_{x \in [0,1]} \|(u_n, t_n) - \tilde{f}_n(x)\| \leq \frac{1}{\sqrt{n}}$. This gives us
    \[
    \|t_n - 1\| \leq \min_{x \in [0,1]} \|(u_n, t_n) - \tilde{f}_n(x)\| \leq \frac{1}{\sqrt{n}} \to 0 \quad \text{as } n \to \infty.
    \]
    \item $|\tilde{\Phi}'(u_n, t_n)(0,1)| \leq \frac{2}{\sqrt{n}} \|(0,1)\| = \frac{2}{\sqrt{n}} \to 0$ as $n \to \infty$.
\end{enumerate}

A direct computation shows that
\begin{equation}\label{D2}
\left| \frac{\partial}{\partial t} \tilde{\Phi}(u_n, t_n) \right| = |\tilde{\Phi}'(u_n, t_n)(0,1)| \to 0 \quad \text{as } n \to \infty.
\end{equation}
Set $v_n := t_n \star u_n$, then by (\ref{D1}) we get $\Phi(v_n) = \Phi(t_n \star u_n) = \tilde{\Phi}(u_n, t_n) \to \tilde{\sigma}_\rho $. By Lemma \ref{L13}, we have $\Phi(v_n) \rightarrow \sigma_\rho$. Also, by (\ref{D2})
\begin{align*}
0 &= \lim_{n \to \infty} \frac{\partial}{\partial t} \tilde{\Phi}(u_n, t_n)\\
&= \lim_{n \to \infty} \left(  a t_n |\nabla u_n|_2^2 + b t_n^{3} |\nabla u_n|_2^4 - \mu \gamma_q t_n^{q \gamma_q -1} |u_n|_q^q - t_n^{2\alpha +5 } \int_{\mathbb{R}^3} (I_\alpha * |u_n|^{\alpha +3}) |u_n|^{\alpha +3} \, dx. \right)\\
&= \lim_{n \to \infty} \frac{1}{t_n} P(t_n \star u_n) \\
&= \lim_{n \to \infty} P(v_n).
\end{align*}

Besides, for $u \in \mathcal{\mathcal{S}}_{\rho , r}$, setting
\[
T_u := \left\{ h \in H_r ^1(\mathbb{R}^3) : \int_{\mathbb{R}^3} u h = 0 \right\},
\]
For any $h_n \in T_{v_n}$, by simple calculate, we have
\begin{align*}
\Phi'(v_n)(h_n) &= t_n^{-\frac{1}{2}}(a+b t_n ^2 |\nabla u_n|_2^2 ) \int_{\mathbb{R}^3} \nabla u_n(y) \cdot \nabla h_n(t_n^{-1} y) \, dy  \\
& \quad -\mu t_n ^{\frac{3(q-3)}{2}} \int_{\mathbb{R}^3} |u_n(y)|^{q-2}u_n(y)h_n(t_n^{-1} y) \, dy\\
&\quad - t_n^{3(\alpha +3)-\frac{9}{2}-\alpha} \int_{\mathbb{R}^3} \int_{\mathbb{R}^3} \frac{A_\alpha |u_n(x)|^{\alpha +3} |u_n(y)|^{\alpha +1} u_n(x) h_n(t_n^{-1} x) }{|x - y|^{3-\alpha}} \, dx \, dy.
\end{align*}
which implies that
\begin{align*}
\Phi'(v_n)(h_n) &= t_n(a+b t_n ^{2} |\nabla u_n|_2^2 ) \int_{\mathbb{R}^3} \nabla v_n(x) \cdot \nabla \hat{h}_n(x) \, dx  \\
& \quad -\mu t_n ^{\frac{3(q-2)}{2}} \int_{\mathbb{R}^3} |u_n(x)|^{q-2}u_n(x)\hat{h}_n(x) \, dx\\
&\quad - t_n^{2(\alpha +3)} \int_{\mathbb{R}^3} \int_{\mathbb{R}^3} \frac{A_\alpha |u_n(x)|^{\alpha +3} |u_n(y)|^{\alpha +1} u_n(x) \hat{h}_n(x) }{|x - y|^{3-\alpha}} \, dx \, dy\\
&= \Phi'(t_n \star u_n)(t_n \star \hat{h}_n) \\
&= \tilde{\Phi}'(u_n, t_n)(\hat{h}_n, 0).
\end{align*}
where $\hat{h}_n(x) = t_n^{-\frac{3}{2}} h_n(t_n^{-1} x)$.

Now, we claim that $(\hat{h}_n, 0) \in \tilde{T}_{(u_n, t_n)}$. In fact
\begin{align*}
(\hat{h}_n, 0) \in \tilde{T}_{(u_n, t_n)} &\Leftrightarrow \langle u_n, \hat{h}_n \rangle=0\\
&\Leftrightarrow \int_{\mathbb{R}^3} u_n t_n^{-\frac{3}{2}}h_n(t_n^{-1}x) \, dx =0\\
&\Leftrightarrow \int_{\mathbb{R}^3} t_n^{\frac{3}{2}} u_n(t_n x) h_n(x) \, dx =0\\
&\Leftrightarrow \langle v_n, h_n \rangle=0\\
&\Leftrightarrow h_n \in T_{v_n}.
\end{align*}
Therefore, by Lemma \ref{L14},
\[
|\Phi'(v_n)(h_n)| = |\tilde{\Phi}'(u_n, t_n)(\hat{h}_n, 0)| \leq \frac{2}{\sqrt{n}} \|(\hat{h}_n, 0)\| \leq \frac{2}{\sqrt{n}} C \|h_n\|^2 \quad \text{for large } n.
\]
Thus,
\[
\sup \{ |\Phi'(v_n)(h)| : h \in T_{v_n}, \|h\| \leq 1 \} \leq \frac{C'}{\sqrt{n}} \to 0 \quad \text{as } n \to \infty.
\]
\end{proof}

Define
\[
m_\rho := \inf_{u \in \mathcal{P}} \Phi(u),
\]
and
\[
m_{\rho, r} := \inf_{u \in \mathcal{P}_r} \Phi(u).
\]
The following results then hold:
\begin{lemma}\label{L16}
For any $\frac{14}{3} < q < 6$, one has
\begin{equation}
m_\rho > 0.
\end{equation}
\end{lemma}
\begin{proof}
For any $u \in \mathcal{P}$, we have
\begin{equation}\label{C1}
\aligned
a |\nabla u|_2^2 + b |\nabla u|_2^4  &= \mu \gamma_q \int_{\mathbb{R}^3} |u|^q \, dx + \int_{\mathbb{R}^3} (I_\alpha * |u|^{\alpha +3}) |u|^{\alpha +3} \, dx\\
&\geq \mu \gamma_q \int_{\mathbb{R}^3} |u|^q \, dx.
\endaligned
\end{equation}
Since $\frac{14}{3}<q<6$, by (\ref{C1}) and Lemma \ref{L10}, we deduce that
\begin{equation}\label{C2}
\aligned
\Phi(u) &=\Phi(u)- \frac{1}{2(\alpha +3)} P(u)\\
&=a \left( \frac{1}{2} - \frac{1}{2(\alpha +3)} \right) |\nabla u|_2^2 + b \left( \frac{1}{4} - \frac{1}{2(\alpha +3)} \right) |\nabla u|_2^4 - \mu \left( \frac{1}{q} - \frac{\gamma_q}{2(\alpha +3)} \right) |u|_q^q .\\
&\geq a \left( \frac{1}{2} - \frac{1}{2(\alpha +3)} \right) |\nabla u|_2^2 + b \left( \frac{1}{4} - \frac{1}{2(\alpha +3)} \right) |\nabla u|_2^4 - \left( \frac{1}{q \gamma_q} - \frac{1}{2(\alpha +3)} \right)(a |\nabla u|_2^2 + b |\nabla u|_2^4)\\
&=a \left( \frac{1}{2} - \frac{1}{q \gamma_q} \right) |\nabla u|_2^2 + b \left( \frac{1}{4} - \frac{1}{q \gamma_q} \right) |\nabla u|_2^4\\
&\geq a \left( \frac{1}{2} - \frac{1}{q \gamma_q} \right) \delta^2 + b \left( \frac{1}{4} - \frac{1}{q \gamma_q} \right) \delta^4 >0.
\endaligned
\end{equation}
Thus, we get $m_\rho > 0$.
\end{proof}

\begin{lemma}\label{L17}
$m_\rho = m_{\rho, r}$.
\end{lemma}
\begin{proof}
For any given $u \in \mathcal{P}$, let $u^*$ denote its symmetric decreasing rearrangement. By the properties of decreasing rearrangement $\Phi(u^*) \leqslant \Phi(u)$ and hence
\[
m_{\rho, r} = \inf_{u \in \mathcal{P}_r} \Phi(u) \geqslant \inf_{u \in \mathcal{P}} \Phi(u) = m_\rho.
\]
On the other hand, for any $\lambda > 0$,
\begin{align*}
&\left| \left\{ x : (t \star u^*)(x) > \lambda \right\} \right| = \left| \left\{ x : t^{\frac{3}{2}} u^*(t x) > \lambda \right\} \right| = \left| \left\{ t^{-1} y : u^*(y) > t^{-\frac{3}{2}} \lambda \right\} \right| \\
&= \left| \left\{ t^{-1} y : u(y) > t^{-\frac{3}{2}} \lambda \right\} \right| = \left| \left\{ x : (t \star u)(x) > \lambda \right\} \right|= \left| \left\{ x : (t \star u)^\ast (x) > \lambda \right\} \right|,
\end{align*}
one can see that
\[
(t \star u)^* = t \star u^*.
\]
Then for any $u \in \mathcal{P}$, we have
\begin{equation}\label{B1}
\Phi(t \star u^*) = \Phi((t \star u)^*) \leqslant \Phi(t \star u) \leqslant \max_{t > 0} \Phi(t \star u) = \Phi(u) \quad \forall u \in \mathcal{P}.
\end{equation}
Moreover, for every $u \in \mathcal{P}$, it holds that $u^* \in \mathcal{S}_{\rho, r}$. According to Lemma \ref{L4}, there exists unique $t_u^* > 0$ such that $t_u^* \star u^* \in \mathcal{P}_r$. Now, using (\ref{B1}) we obtain $\Phi(t_u^* \ast u^*) \leqslant \Phi(u)$ and
\[
m_{\rho, r} = \inf_{u \in \mathcal{P}_r} \Phi(u) \leq \inf_{u \in \mathcal{P}} \Phi(u) = m_\rho.
\]
Hence, we have $m_{\rho, r} = m_\rho.$
\end{proof}

\begin{lemma}\label{L18}
$m_{\rho, r} = \sigma_\rho > 0$.
\end{lemma}
\begin{proof}
First, we assert that if $\Phi(u) \leqslant 0$ for some $u \in \mathcal{S}_{\rho, r}$, then $P(u) < 0$. Consider $u \in \mathcal{S}{\rho, r}$ satisfying $\Phi(u) \leqslant 0$. By Lemma \ref{L4}, there exists a unique $t_u \in \mathbb{R}^+$ such that $t_u \star u \in \mathcal{P}_r$. Note that $E_u(0) = 0$, with $E_u$ strictly increasing on $(0, t_u)$ and strictly decreasing on $(t_u, \infty)$, which implies $E_u(t_u) > 0$. Since $E_u(1) = \Phi(u) \leqslant 0$, it follows that $t_u < 1$. Thus, according to Lemma \ref{L2}, we obtain $P(u) = E'_u(1) < 0$.

Next, we show that $\sigma_\rho = m_{\rho, r}$.
For $u \in \mathcal{P}_r \subset \mathcal{S}_{\rho, r}$, let $t^- \in (0,1)$ and $t^+ > 1$ be such that $t^- \star u \in \overline{A_k}$ and $\Phi(t^+ \star u) < 0$. Define the path $\eta_u : [0,1] \rightarrow \mathcal{S}_{\rho, r}$ by
\[
\eta_u(z) = ((1 - z)t^- + zt^+) \ast u.
\]
Since $\eta_u \in \Gamma_\rho$, it follows that $\sigma_\rho \leqslant \max_{z \in [0,1]} \Phi(\eta_u(z))$. Moreover, according to Remark \ref{L5},
\[
\Phi(u) = E_u(1) \geqslant E_u(t) = \Phi(t \star u) \quad \text{for all } t \in [t^-, t^+].
\]
Therefore
\[
\Phi(u) = \max_{z \in [0,1]} \Phi(\eta_u(z)) \geqslant \sigma_\rho,
\]
since $u \in \mathcal{P}_r$ is arbitrary, we get
\[
m_{\rho, r} \geq \sigma_\rho .
\]
Now for any $\tilde{\eta} = (\tilde{\eta}_1, \tilde{\eta}_2) \in \tilde{\Gamma}_\rho$, let us consider the function
\[
\tilde{P}(z) : z \in [0,1] \mapsto P(\tilde{\eta}_2(z) \star \tilde{\eta}_1(z)),
\]
By Lemma \ref{L12}, $\tilde{P}(0) = P(\tilde{\eta}_1(0)) > 0$ and since $\tilde{\eta}_1(1) \in \Phi^0$, we get $\Phi(\tilde{\eta}_1(1)) \leqslant 0$. Thus by the first claim, $\tilde{P}(1) = P(\tilde{\eta}_1(1)) < 0$. Hence, there exists $\tilde{z} \in (0,1)$ such that $\tilde{P}(\tilde{z}) = 0$, which implies that $\tilde{\eta}_2(\tilde{z}) \star \tilde{\eta}_1(\tilde{z}) \in \mathcal{P}_r$. Therefore, we get
\[
\max_{z \in [0,1]} \Phi(\tilde{\eta}_2(z) \star \tilde{\eta}_1(z)) \geqslant \Phi(\tilde{\eta}_2(\tilde{z}) \star \tilde{\eta}_1(\tilde{z})) \geqslant \inf_{u \in P_r(\tau)} \Phi(u) = m_{\rho, r}.
\]
Since $\tilde{\eta} \in \tilde{\Gamma}_\rho$ is arbitrary, we get $\sigma_\rho \geqslant m_{\rho, r}$. Therefore, we have $m_{\rho, r} = \sigma_\rho > 0$.
\end{proof}

For ease of calculation, in the following proof, we denote
\[
\mathcal{A}(u)= \int_{\mathbb{R}^3} (I_\alpha * |u|^{\alpha +3}) |u|^{\alpha +3} \, dx.
\]
\begin{lemma}\label{L19}
There exists $ \kappa, \mu^* > 0$ such that for $0 < \rho < \kappa$ and $\mu > \mu^*$
\[
\sigma_\rho < \left( \frac{\alpha +2}{2(\alpha +3)} \right) (aS_\alpha)^{\frac{\alpha +3}{\alpha +2}}.
\]
\end{lemma}
\begin{proof}
By \cite{GHPS} and \cite{SM}, we know that $S_\alpha$ is attained by
\[
U_\epsilon(x) = \frac{(3\epsilon^2)^{\frac{1}{4}}}{(\epsilon^2 + |x|^2)^{\frac{1}{2}}}, \quad \forall \varepsilon >0.
\]
Take a radially decreasing cut-off function $\psi \in C_c^\infty(\mathbb{R}^3, [0, 1])$ satisfy $\psi(x) = 1$ in $B_1(0)$, $\psi(x) = 0$ in $ \mathbb{R}^3 \setminus B_2(0)$, and let $u_\epsilon(x) := \psi(x) U_\epsilon(x)$, $v_\epsilon(x) := \frac{\rho \cdot u_\epsilon(x)}{|u_\epsilon|_2} \in \mathcal{S}_\rho$. We have
\[
|\nabla v_\epsilon|_2^2 = \frac{\rho^2}{|u_\epsilon|_2^2}|\nabla u_\epsilon|_2^2,
\]
\[
|v_\epsilon|_q^q=\frac{\rho^q}{|u_\epsilon|_2^q}|u_\epsilon|_q^q,
\]
and
\[
\mathcal{A}(v_\epsilon)=\frac{\rho^{2(\alpha+3)}}{|u_\epsilon|_2^{2(\alpha+3)}}A(u_\epsilon).
\]
As shown in Lemma 4.7 of \cite{SM}:
\begin{equation}\label{A1}
|\nabla u_\epsilon|_2^2 = S^{\frac{3}{2}} + O(\epsilon),
\end{equation}
\begin{equation}\label{A2}
|u_\epsilon|_2^2 = C_1 \epsilon + O(\epsilon^2),
\end{equation}
\begin{equation}\label{A3}
\mathcal{A}(u_\epsilon) \geqslant (\mathcal{A}_\alpha C_\alpha)^{\frac{3}{2}} S_\alpha^{\frac{\alpha +3}{2}} - O(\epsilon^{\frac{\alpha +3}{2}}),
\end{equation}
and
\begin{equation}\label{A4}
|u_\epsilon|_p^p = \epsilon^{3 - \frac{p}{2}} (C_2 + O(\epsilon^{p-3})),
\end{equation}
for some positive constant $C_1$ and $C_2$. Since $v_\epsilon \in \mathcal{S}_\rho$, it follows from Lemma \ref{L4} there exists unique $t_\epsilon > 0$ satisfying $E'_{v_\epsilon}(t_\epsilon) = 0$. Define
\[
\mathcal{B}_\epsilon(u) := a|\nabla u|_2^2 + b t_\epsilon^{2} |\nabla u|_2^4,
\]
from this, we obtain
\[
t_\epsilon^{2\alpha +5} \mathcal{A}(v_\epsilon) = t_\epsilon \mathcal{B}_\epsilon(v_\epsilon) - \mu \gamma_q t_\epsilon^{\frac{3(q-2)}{2}-1} |v_\epsilon|_q^q.
\]
Which implies
\begin{equation}\label{A5}
t_\epsilon \leqslant \left( \frac{\mathcal{B}_\epsilon(v_\epsilon)}{\mathcal{A}(v_\epsilon)} \right)^{\frac{1}{2 \alpha +4}}.
\end{equation}
Moreover, since $E'_{v_\epsilon}(t_\epsilon) = 0$, it follows that
\begin{align*}
t_\epsilon^{2\alpha +4} &= \frac{\mathcal{B}_\epsilon(v_\epsilon)}{\mathcal{A}(v_\epsilon)} - \frac{\mu \gamma_q t_\epsilon^{\frac{3(q-2)}{2}-2} |v_\epsilon|_q^q}{\mathcal{A}(v_\epsilon)}\\
&=\frac{\mathcal{B}_\epsilon(v_\epsilon)}{2\mathcal{A}(v_\epsilon)} + \left( \frac{\mathcal{B}_\epsilon(v_\epsilon)}{2\mathcal{A}(v_\epsilon)} - \frac{\mu \gamma_q t_\epsilon^{\frac{3(q-2)}{2}-2} |v_\epsilon|_q^q}{\mathcal{A}(v_\epsilon)} \right) .
\end{align*}
Using the definition of $\mathcal{B}_\epsilon(u)$ and $v_\epsilon(x)$, it follows that
\begin{align*}
\mathcal{B}_\epsilon(v_\epsilon) &= \frac{\rho^2}{|u_\epsilon|_2^2}a|\nabla u_\epsilon|_2^2 +bt_\epsilon^2\frac{\rho^4}{|u_\epsilon|_2^4}|\nabla u_\epsilon|_2^4\\
&=\frac{\rho^2}{|u_\epsilon|_2^2}(a|\nabla u_\epsilon|_2^2 + bt_\epsilon^2\frac{\rho^2}{|u_\epsilon|_2^2}|\nabla u_\epsilon|_2^4).
\end{align*}

Define
\[
\mathcal{D}_\epsilon(u_\epsilon) := a|\nabla u_\epsilon|_2^2 + bt_\epsilon^2\frac{\rho^2}{|u_\epsilon|_2^2}|\nabla u_\epsilon|_2^4.
\]
We claim that for $0 < \epsilon < \epsilon_0$
\begin{equation}\label{A6}
\frac{\mathcal{B}_\epsilon(v_\epsilon)}{2\mathcal{A}(v_\epsilon)} - \frac{\mu \gamma_q t_\epsilon^{\frac{3(q-2)}{2}-2} |v_\epsilon|_q^q}{\mathcal{A}(v_\epsilon)} > 0,
\end{equation}
which implies
\begin{equation*}
t_\epsilon^{2\alpha +4} > \frac{\mathcal{B}_\epsilon(v_\epsilon)}{2\mathcal{A}(v_\epsilon)} > \frac{a|\nabla v_\epsilon|_2^2}{2\mathcal{A}(v_\epsilon)} = \frac{a|u_\epsilon|_2^{2\alpha +4} |\nabla u_\epsilon|_2^2}{2 \rho^{2\alpha +4} \mathcal{A}(u_\epsilon)} \geqslant C_{a, \rho} S_\alpha^{\alpha +3} |u_\epsilon|_2^{2\alpha +4},
\end{equation*}
that is,
\begin{equation}\label{A7}
t_\epsilon > C_{a, \rho, \alpha} |u_\epsilon|_2 \quad \text{for } 0 < \epsilon < \epsilon_0.
\end{equation}
From (\ref{A5}) and the definition of $v_\epsilon$, it follows that
\begin{align*}
\mathcal{B}_\epsilon(v_\epsilon) &- 2 \mu \gamma_q t_\epsilon^{\frac{3(q-2)}{2}-2} |v_\epsilon|_q^q \geqslant \frac{\rho^2 \mathcal{D}_\epsilon(u_\epsilon)}{|u_\epsilon|_2^2} - \frac{2\mu \gamma_q |u_\epsilon|_q^q \rho^q}{|u_\epsilon|_2^q} \left( \frac{|u_\epsilon|_2^{2\alpha +4} \mathcal{D}_\epsilon(u_\epsilon)}{\rho^{2\alpha +4} \mathcal{A}(u_\epsilon)} \right)^{\frac{3(q-2)-4}{4(\alpha +2)}}\\
&= \frac{\rho^2}{|u_\epsilon|_2^2} \mathcal{D}_\epsilon(u_\epsilon)^{\frac{3(q-2)-4}{4(\alpha +2)}} \left( \mathcal{D}_\epsilon(u_\epsilon)^{\frac{4(\alpha +3)-3(q-2)}{4(\alpha +2)}} - \frac{2\mu \gamma_q \rho^{q-\frac{3(q-2)}{2}}}{\mathcal{A}(u_\epsilon)^{\frac{3(q-2)-4}{4(\alpha +2)}}} \frac{|u_\epsilon|_q^q}{|u_\epsilon|_2^{q-\frac{3(q-2)}{2}}}\right).
\end{align*}
Thus, (\ref{A6}) holds provided that
\begin{equation}\label{A8}
2\mu \gamma_q \rho^{q-\frac{3(q-2)}{2}} < \frac{\mathcal{A}(u_\epsilon)^{\frac{3(q-2)-4}{4(\alpha +2)}} \mathcal{D}_\epsilon(u_\epsilon)^{\frac{4(\alpha +3)-3(q-2)}{4(\alpha +2)}}|u_\epsilon|_2^{q-\frac{3(q-2)}{2}}}{|u_\epsilon|_q^q}
\quad \text{for } 0 < \epsilon < \epsilon_0.
\end{equation}
Using estimates (\ref{A1}), (\ref{A3}), and (\ref{A4}) for $0 < \epsilon < \epsilon_0$, we obtain
\begin{align*}
\frac{\mathcal{A}(u_\epsilon)^{\frac{3(q-2)-4}{4(\alpha +2)}} \mathcal{D}_\epsilon(u_\epsilon)^{\frac{4 (\alpha +3) - 3(q-2)}{4(\alpha +2)}}|u_\epsilon|_2^{q-\frac{3(q-2)}{2}}}{ |u_\epsilon|_q^q} &\geqslant C \frac{\mathcal{D}_\epsilon(u_\epsilon)^{\frac{2 q - 3(q-2)}{4(\alpha +2)}}\mathcal{A}(u_\epsilon)^{\frac{3(q-2)-4}{4(\alpha +2)}} |u_\epsilon|_2^{q-\frac{3(q-2)}{2}}}{|u_\epsilon|_q^q} \\
&\geqslant C(\alpha, q) \left(\frac{a|\nabla u_\epsilon|_2^2}{\mathcal{A}(u_\epsilon)^{\frac{1}{\alpha + 3}}}\right)^{\frac{2 q - 3(q-2)}{4(\alpha +2)}}\\
&\geqslant C(\alpha, q, a) S_\alpha ^{\frac{2 q - 3(q-2)}{4(\alpha +2)}}.
\end{align*}
Hence, for
\begin{equation}\label{A9}
\rho < \kappa := \left( \frac{C(\alpha, q, a)}{2\mu \gamma_q} \right)^{\frac{2}{6-q}} S_\alpha^{\frac{1}{2\alpha +4}},
\end{equation}
we get
\[
2\mu \gamma_q \rho^{q - \frac{3(q-2)}{2}} < C(\alpha, q, a) S_\alpha^{\frac{2 q - 3(q-2)}{4(\alpha +2)}} \leqslant \frac{\mathcal{A}(u_\epsilon)^{\frac{3(q-2)-4}{4(\alpha +2)}} \mathcal{D}_\epsilon(u_\epsilon)^{\frac{4(\alpha +3) - 3(q-2)}{4(\alpha +2)}}|u_\epsilon|_2^{q-\frac{3(q-2)}{2}}}{ |u_\epsilon|_q^q} .
\]

Define
\[
I_{v_\epsilon}(t) := \frac{at^2}{2} |\nabla v_\epsilon|_2^2 - \frac{t^{2(\alpha +3)}}{2(\alpha +3)} \mathcal{A}(v_\epsilon).
\]
Obviously, $I_{v_\epsilon}$ has its maximum at $t_0 = \left( \frac{a|\nabla v_\epsilon|_2^2}{\mathcal{A}(v_\epsilon)} \right)^{\frac{1}{2\alpha +4}}$. By (\ref{A1}) and (\ref{A3}) we have
\begin{align*}
I_{v_\epsilon}(t_0) &= \frac{\alpha +2}{2 (\alpha +3)} \left( \frac{a|\nabla u_\epsilon|_2^2}{\mathcal{A}(u_\epsilon)^{\frac{1}{\alpha +3}}} \right)^{\frac{\alpha +3}{\alpha +2}}\\
& \leqslant a^{\frac{\alpha +3}{\alpha +2}}\left( \frac{\alpha +2}{2 (\alpha +3)} \right) \left( \frac{S^{\frac{3}{2}} + O(\epsilon)}{((\mathcal{A}_\alpha C_\alpha)^{\frac{3}{2}} S_\alpha^{\frac{\alpha +3}{2}} - O(\epsilon^{\frac{\alpha +3}{2}}))^{\frac{1}{\alpha +3}}} \right)^{\frac{\alpha +3}{\alpha +2}}\\
&=a^{\frac{\alpha +3}{\alpha +2}} \left( \frac{\alpha +2}{2(\alpha +3)} \right) \left( \frac{(\mathcal{A}_\alpha C_\alpha)^{\frac{3}{2(\alpha +3)}} S_\alpha^{\frac{3}{2}} + O(\epsilon)}{(\mathcal{A}_\alpha C_\alpha)^{\frac{3}{2(\alpha +3)}} S_\alpha^{\frac{\alpha +3}{2(\alpha +3)}} (1 - O(\epsilon^{\frac{\alpha +3}{2}}))^{\frac{1}{\alpha +3}}} \right)^{\frac{\alpha +3}{\alpha +2}} \\
&\leqslant \left( \frac{\alpha +2}{2 (\alpha +3)} \right) \left(aS_\alpha\right)^{\frac{\alpha +3}{\alpha +2}} + O(\epsilon).
\end{align*}
Therefore
\[
\sup_{t > 0} I_{v_\epsilon}(t) =a^{\frac{\alpha +3}{\alpha +2}} \left( \frac{\alpha +2}{2 (\alpha +3)} \right) S_\alpha^{\frac{\alpha +3}{\alpha +2}} + O(\epsilon).
\]
Either $t_\epsilon < 1$ or by (\ref{A5}),
\[
t_\epsilon \leqslant \left( \frac{\mathcal{B}_\epsilon (v_\epsilon)}{\mathcal{A}(v_\epsilon)} \right)^{\frac{1}{2\alpha +4}} = \left( \frac{a|\nabla v_\epsilon|_2^2 + b t_\epsilon^{2} |\nabla v_\epsilon|_2^4}{\mathcal{A}(v_\epsilon)} \right)^{\frac{1}{2\alpha +4}} \leq t_\epsilon ^{\frac{1}{\alpha +2}}\left( \frac{a|\nabla v_\epsilon|_2^2 + b |\nabla v_\epsilon|_2^4}{\mathcal{A}(v_\epsilon)} \right)^{\frac{1}{2\alpha +4}}.
\]
So, we have $t_\epsilon \leq \vartheta_1$ for $0 < \epsilon < \epsilon_0$.
The existence of this $\vartheta_1$ can be deduced using (\ref{A1}), (\ref{A2}), (\ref{A3}) and (\ref{A4}). Taking $\vartheta = \max\{1, \vartheta_1\}$, we get $t_\epsilon \leqslant \vartheta$. Thus we get
\begin{align*}
\sup_{t > 0} E_{v_\epsilon}(t) &= E_{v_\epsilon}(t_\epsilon) = I_{v_\epsilon}(t_\epsilon) + \frac{bt_\epsilon^{4}}{4} |\nabla v_\epsilon|_2^4 - \frac{\mu}{q}t_\epsilon^{\frac{3(q-2)}{2}} |v_\epsilon|_q^q\\
&\leqslant \left( \frac{\alpha +2}{2(\alpha +3)} \right) \left(aS_\alpha\right)^{\frac{\alpha +3}{\alpha +2}} + O(\epsilon) + \frac{b\vartheta^{4}}{4} |\nabla v_\epsilon|_2^4 - \frac{\mu}{q}t_\epsilon^{\frac{3(q-2)}{2}} |v_\epsilon|_q^q.
\end{align*}
By the definition of $v_\epsilon$, we have
\begin{equation*}
\frac{b\vartheta^{4}}{4} |\nabla v_\epsilon|_2^4 - \frac{\mu}{q}t_\epsilon^{\frac{3(q-2)}{2}} |v_\epsilon|_q^q = \frac{b\vartheta^{4} \rho^4 |\nabla u_\epsilon|_2^4}{4| u_\epsilon|_2^4} - \frac{\mu}{q} \frac{ t_\epsilon^{\frac{3(q-2)}{2}} \rho^q |u_\epsilon|_q^q}{|u_\epsilon|_2^{q}}.
\end{equation*}

Set $Y_\epsilon := \frac{|\nabla u_\epsilon|_2^4}{|u_\epsilon|_2^4}$ and $Y^{'}_\epsilon := \frac{t_\epsilon^{\frac{3(q-2)}{2}} |u_\epsilon|_q^q}{|u_\epsilon|_2^q}$.
By Lemma \ref{L16} we have
\[
\Phi(u) \geqslant m_\rho = m_{\rho, r} \geqslant a \left( \frac{1}{2} - \frac{1}{q \gamma_q} \right) \delta^2 + b \left( \frac{1}{4} - \frac{1}{q \gamma_q} \right) \delta^4 \quad \text{for all } u \in \mathcal{P}.
\]
Let $0 < l < a \left( \frac{1}{2} - \frac{1}{q \gamma_q} \right) \delta^2 + b \left( \frac{1}{4} - \frac{1}{q \gamma_q} \right) \delta^4$, then $l < \Phi(t_\epsilon \star v_\epsilon) = \max_{t \geqslant 0} \Phi(t \star u_\epsilon) = E_{v_\epsilon}(t_\epsilon)$ for every $\epsilon > 0$. Also, since $E_{v_\epsilon}(0) = 0$ for all $\epsilon > 0$, then there must exist $t_0 > 0$ small enough such that $t_\epsilon \geqslant t_0$ for all $\epsilon > 0$. By (\ref{A1}), (\ref{A2}), (\ref{A4}) and (\ref{A7}),
\[
\frac{\epsilon t_\epsilon^4 Y_\epsilon}{Y_\epsilon^{'}} \leqslant |\nabla u_\epsilon|_2^4 \epsilon^{\frac{q}{4} - \frac{1}{2}} \to 0 \quad \text{as } \epsilon \to 0.
\]
Therefore, by (\ref{A1}) we obtain that $\frac{t_\epsilon^4 Y_\epsilon}{Y_\epsilon^{'}} \leqslant S^3$ and hence $\frac{Y_\epsilon}{Y_\epsilon^{'}} \leqslant \frac{S^3}{t_0^4}$. From the above, it follows that
\begin{equation}
\frac{b\vartheta^{4}}{4} |\nabla v_\epsilon|_2^4 - \frac{\mu}{q}t_\epsilon^{\frac{3(q-2)}{2}} |v_\epsilon|_q^q <0 \quad \text{for } \epsilon \to 0,
\end{equation}
if
\begin{equation}\label{A10}
\mu >\frac{bqS^3\vartheta^4}{t_0^4}:=\mu ^*.
\end{equation}
This gives us
\[
\sup_{t > 0} E_{v_\epsilon}(t) < \left( \frac{\alpha +2}{2(\alpha +3)} \right) \left(aS_\alpha \right)^{\frac{\alpha +3}{\alpha +2}}
\]
as $\epsilon \to 0$. Define the path $\eta_{u_\epsilon}(z) = ((1 - z) t^- + z t^+) \star v_\epsilon$ for $z \in [0,1]$, where $t^-, t^+$ are as in the proof of Lemma \ref{L18}. Since $\eta_{u_\epsilon} \in \Gamma_\rho$, we have
\[
\sigma_\rho \leqslant \sup_{z \in [0,1]} E_{v_\epsilon}(z) \leqslant \sup_{t > 0} E_{v_\epsilon}(t) < \left( \frac{\alpha +2}{2 (\alpha +3)} \right) \left(aS_\alpha \right)^{\frac{\alpha +3}{\alpha +2}}.
\]
\end{proof}

\section{Proof of Theorem 1.1}
Now we prove the existence of normalized solution for Eq.(\ref{p})-(\ref{q}) with $\frac{14}{3} < q <6$.

{\bf Proof of Theorem 1.1}\ Let $\{u_n\}$ be the sequence in $\mathcal{P}_r$ satisfying
\[
\lim_{n \to \infty} \Phi(u_n) = m_{\rho,r} = \inf_{u \in \mathcal{P}_r} \Phi(u).
\]
By (\ref{C2}), $\Phi$ is coercive on $\mathcal{P}_r$, so $\{u_n\}$ is bounded and possesses a weakly convergent subsequence (still denoted by $\{u_n\}$) in $H^1_r(\mathbb{R}^3)$. Denote by $u_\rho \in H_r^1(\mathbb{R}^3)$ the weak limit of $\{u_n\}$. Since $\{u_n\}$ is a minimizing sequence for $\Phi$ on $\mathcal{P}_r$, the Ekeland variational principle \cite{J} implies the existence of sequences ${\lambda_n}$ and ${\zeta_n}$ such that
\[
\Phi'(u_n) - \lambda_n u_n + \zeta_n P'(u_n) \to 0 \quad \text{as } n \to \infty.
\]
Define $J_n(u) := \Phi(u) - \lambda_n |u|2^2 + \zeta_n P(u)$; then $J_n'(u_n) \to 0$. We will show that $\zeta_n \to 0$ and $\lambda_n \to \lambda\rho$ for some $\lambda_\rho$, so that $(\lambda_\rho, u_\rho)$ yields the desired result.

\textbf{Step 1:} $\zeta_n \rightarrow 0$ .

Since $u_n \in \mathcal{P}_r$ for every $n \in \mathbb{N}$, we have:
\begin{align*}
0 &= \lim_{n \to \infty} \left. \frac{\partial}{\partial t} J_n(t \star u_n) \right|_{t=1} \\
&= \lim_{n \to \infty} \frac{\partial}{\partial t} \left (\Phi(t \star u_n) - \lambda_n |u_n|_2^2 + \zeta_n P(t \star u_n)\right )_{t=1} \\
&= \lim_{n \to \infty} (E'_{u_n}(1) + \zeta_n (E''_{u_n}(1) + E'_{u_n}(1))) \\
&= \lim_{n \to \infty} \zeta_n E''_{u_n}(1),
\end{align*}
also Lemma \ref{L3} gives us $|E''_{u_n}(1)| = -E''_{u_n}(1) \geqslant \left(q\gamma _q-2\right)a\delta ^2+\left(q\gamma _q-4\right)b\delta ^4$, then
\[
|\zeta_n \delta^2| \leqslant \frac{ |\zeta_n E''_{u_n}(1)|}{a(q\gamma_q - 2)} \to 0 \quad \text{as } n \to \infty,
\]
as $\delta > 0$ for all $n \in \mathbb{N}$, we get
\[
\zeta_n \to 0 \quad \text{as } n \to \infty.
\]

\textbf{Step 2:} $u_\rho \neq 0$.

Let if possible $u_\rho = 0$, taking $l \in \mathbb{R}$ such that $ |\nabla u_n|_2^2 \rightarrow l$. Since the embedding $H^1_r(\mathbb{R}^3) \hookrightarrow L^p(\mathbb{R}^3)$ is compact for all $p \in (2, 6)$, $\{u_n\} \to u_\rho$ in $L^p(\mathbb{R}^3)$ for every $p \in (2, 6)$. Further, since $P(u_n) = 0$ for every $n \in \mathbb{N}$, we get
\begin{equation}\label{E31}
\aligned
\lim_{n \to \infty} \mathcal{A}(u_n) &= \lim_{n \to \infty} \left( a|\nabla u_n|_2^2 +b|\nabla u_n|_2^4 - \mu\gamma_q |u_n|_q^q - P(u_n) \right) \\
&= \lim_{n \to \infty} \left( a|\nabla u_n|_2^2 +b|\nabla u_n|_2^4 \right) \\
&= al + bl^2.
\endaligned
\end{equation}
Now, since
\[
S_\alpha \leqslant \frac{|\nabla u_n|_2^2}{\left( \int_{\mathbb{R}^3} (I_\alpha * |u_n|^{3+\alpha}) |u_n|^{3+\alpha} \, dx \right)^{\frac{1}{\alpha +3}}} \quad \text{for every } n \in \mathbb{N},
\]
then by (\ref{E31}) we get
\[
S_\alpha \leqslant \frac{l}{\left(al+bl^2\right)^{\frac{1}{\alpha +3}}}\leq  \frac{l}{\left(al\right)^{\frac{1}{\alpha +3}}} \Rightarrow l (aS_\alpha^{\alpha +3} - l^{\alpha +2}) \leqslant 0,
\]
thus either $l = 0$ or $l \geqslant a^{\frac{1}{\alpha +2}}S_\alpha^{\frac{\alpha +3}{\alpha +2}}$. If $l \geqslant a^{\frac{1}{\alpha +2}}S_\alpha^{\frac{\alpha +3}{\alpha +2}}$, then
\begin{align*}
\sigma_\rho = m_{\rho, r} &= \lim_{n \to \infty} \Phi(u_n) = \lim_{n \to \infty} \left( \Phi(u_n) - \frac{1}{2(\alpha +3)} P(u_n) \right)\\
&= \lim_{n \to \infty} \left( \frac{a}{2} \left( 1 - \frac{1}{\alpha +3} \right) |\nabla u_n|_2^2 + b\left(\frac{1}{4}  - \frac{1}{2(\alpha +3)}\right) |\nabla u_n|_2^4 \right)\\
&\geqslant  \frac{(\alpha +2)a}{2(\alpha +3)} l \geqslant \frac{(\alpha +2)}{2(\alpha +3)} \left(aS_\alpha \right)^{\frac{\alpha +3}{\alpha +2}},
\end{align*}
this contradicts Lemma \ref{L19}. Therefore $l = 0$ and hence $\sigma_\rho = \lim_{n \to \infty} \Phi(u_n) = 0$, but since $\sigma_\rho > 0$, we must have $u_\rho \neq 0$.

\textbf{Step 3:} ${\lambda_n}$ is convergent.

Note that $\{u_n\}$ is a Palais-Smale sequence of $\Phi$, by the boundedness of $\{u_n\}$ and Lagrange multiplier rule, there exists a sequence $\lambda_n \in \mathbb{R}$ such that
\begin{equation}\label{E32}
\langle \Phi'(u_n), \varphi \rangle - \lambda_n \langle u_n, \varphi \rangle = o_n(1), \ \forall \varphi \in H_r^1(\mathbb{R}^3).
\end{equation}
In particular, it follows from $\varphi = u_n$ that
\begin{equation}\label{E33}
a |\nabla u_n|_2^2 + b |\nabla u_n|_2^4 - \mu |u_n|^q_q  - \int_{\mathbb{R}^3} (I_\alpha * |u_n|^{3+\alpha}) |u_n|^{\alpha +3} \, dx - \lambda_n |u_n|^2_2 = o_n(1).
\end{equation}
Further, for $P(u_n) \to 0$, we have
\begin{equation}\label{E34}
P(u_n) = a |\nabla u_n|_2^2 + b |\nabla u_n|_2^4 - \mu \gamma_q |u_n|^q_q - \int_{\mathbb{R}^3} (I_\alpha * |u_n|^{\alpha+3}) |u_n|^{\alpha+3} dx = o_n(1).
\end{equation}
Combining this with (\ref{E33}) and (\ref{E34}), we deduce that
\begin{equation}\label{E35}
\lambda_n |u_n|^2_2 = \mu (\gamma_q - 1) |u_n|^q_q + o_n(1).
\end{equation}
It follows that $\lambda_n\}$ is bounded from $H_r^1(\mathbb{R}^3) \hookrightarrow L^p(\mathbb{R}^3)$ and (\ref{E35}), so up to a subsequence, we may assume that $\lambda_n \to \lambda_\rho \in \mathbb{R}$. Moreover,  since $\{u_n\} $ is bounded, up to a subsequence, there exists a function $u_\rho \in H_r^1(\mathbb{R}^3)$ such that
\begin{equation}\label{E36}
\aligned
& u_n \rightharpoonup u_\rho, \text{ in } H_r^1(\mathbb{R}^3), \\
&u_n \to u_\rho, \text{ in } L^p(\mathbb{R}^3), \ \ \forall p \in (2, 6),\\
&u_n \to u_\rho, \text{ a.e. in } \mathbb{R}^3.
\endaligned
\end{equation}
So we get
\begin{equation}\label{E37}
\lambda_\rho \rho^2 = \mu (\gamma_q - 1) |u_\rho|^q_q.
\end{equation}
Since $u_\rho \neq 0$, we have $\lambda_\rho < 0$.

\textbf{Step 4:} We prove $u_\rho \in S_\rho$.

We show $u_n \to u_\rho$ in $L^2(\mathbb{R}^3)$. In fact, since $u_n \rightharpoonup u_\rho \neq 0$, we may assume that
\begin{equation}\label{E38}
\tilde{B} = \lim_{n \to \infty}  |\nabla u_n|^2_2 \geq |\nabla u_\rho|^2_2 > 0.
\end{equation}
Then, $u_n \rightharpoonup u_\rho$ in $H_r^1(\mathbb{R}^3)$, $H_r^1(\mathbb{R}^3) \hookrightarrow L^q(\mathbb{R}^3)$, (\ref{E32}) and Lemma \ref{L7} imply that
\begin{equation}\label{E39}
(a + \tilde{B}b) \int_{\mathbb{R}^3} \nabla u_\rho \cdot \nabla \varphi - \int_{\mathbb{R}^N} |u_\rho|^{p-2} u_\rho \varphi - \int_{\mathbb{R}^3} (I_\alpha * |u_\rho|^{\alpha +3}) |u_\rho|^{\alpha +1} u_\rho \varphi - \lambda_\rho \int_{\mathbb{R}^3} u_\rho \varphi = 0
\end{equation}
for any $\varphi \in H^1(\mathbb{R}^3)$. That is,
\begin{equation}\label{E310}
-(a + \tilde{B}b) \Delta u_\rho = \lambda_\rho u_\rho + \mu  |u_\rho|^{p-2} u_\rho + (I_\alpha * |u_\rho|^{\alpha +3}) |u_\rho|^{\alpha +1} u_\rho.
\end{equation}
By Lemma \ref{L1}, we get
\begin{equation}\label{E311}
(a + \tilde{B}b) |\nabla u_\rho|^2_2 - \mu \gamma_q  |u_\rho|^q _q - \int_{\mathbb{R}^3} (I_\alpha * |u_\rho|^{\alpha +3}) |u_\rho|^{\alpha +1} dx = 0.
\end{equation}
In particular, we take $\varphi = u_\rho$ in (\ref{E39}), then by (\ref{E311}), we have
\[
\lambda_\rho  |u_\rho|^2_2  = \mu (\gamma_q - 1)  |u_\rho|^q _q.
\]
Together with (\ref{E37}), we can see that
\[
\lambda_\rho \left(\rho^2 - |u_\rho|^2_2 \right) = 0.
\]
So in view of $\lambda_\rho < 0$, we have $|u_\rho|_2^2 = \rho^2$, namely $u_\rho \in \mathcal{S}_\rho$.

\textbf{Step 5:} We verify that $u_n \to u_\rho$ in $H_r^1(\mathbb{R}^3)$.

Denote $\tilde{v}_n = u_n - u_\rho$, by Brezis-Lieb lemma of the Hardy-Littlewood-Sobolev upper critical exponent type(\cite{GHPS}, Lemma 3.1) and (\ref{E38}), we have
\begin{equation}\label{E312}
|\nabla u_n|_2^2 = |\nabla \tilde{v}_n|_2^2 + |\nabla u_\rho|_2^2 + o_n(1), \quad |\nabla u_n|_2^2 = \tilde{B} + o_n(1),
\end{equation}
and
\begin{equation}\label{E313}
\aligned
\int_{\mathbb{R}^3} (I_\alpha * |u_n|^{\alpha +3}) |u_n|^{\alpha +3} \, dx =&\int_{\mathbb{R}^3} (I_\alpha * |\tilde{v}_n|^{\alpha +3}) |\tilde{v}_n|^{\alpha +3} \, dx \\
&+\int_{\mathbb{R}^3} (I_\alpha * |u_\rho|^{\alpha +3}) |u_\rho|^{\alpha +3} \, dx + o_n(1).
\endaligned
\end{equation}
Using (\ref{E38}) again, it is easy to see that
\begin{equation}\label{E314}
\aligned
|\nabla u_n|_2^4 &= (\tilde{B} + o_n(1)) (|\nabla u_\rho|_2^2 + |\nabla \tilde{v}_n|_2^2 + o_n(1)) \\
&\geq |\nabla u_\rho|_2^4 + \tilde{B}|\nabla \tilde{v}_n|_2^2 + o_n(1).
\endaligned
\end{equation}
Combining (\ref{E34}) with (\ref{E310}), (\ref{E311}), then by (\ref{E313}), we get
\begin{equation}\label{E315}
(a + \tilde{B}b) |\nabla \tilde{v}_n|_2^2 - \int_{\mathbb{R}^3} (I_\alpha * |\tilde{v}_n|^{\alpha +3}) |\tilde{v}_n|^{\alpha +3} \, dx = o_n(1).
\end{equation}
By virtue of (\ref{E36}), we get
\begin{equation}\label{E316}
|u_n|_q \to |u_\rho|_q, \ 2 < q < 6.
\end{equation}
Therefore, it follows from (\ref{E312}), (\ref{E313}), (\ref{E314}), (\ref{E315}) and (\ref{E316}) that
\begin{align*}
m_\rho + o_n(1) &= \Phi(u_n)\\
&\geq \Phi(u_\rho) + \frac{a}{2} |\nabla \tilde{v}_n|_2^2 + \frac{b \tilde{B}}{4} |\nabla \tilde{v}_n|_2^2 - \frac{1}{2(\alpha +3)} \int_{\mathbb{R}^3} (I_\alpha * |\tilde{v}_n|^{\alpha +3}) |\tilde{v}_n|^{\alpha +3} \, dx + o_n(1)\\
& \geq m_\rho + \left( \frac{a}{2} - \frac{a}{2(\alpha +3)} + \frac{\tilde{B} b}{4} - \frac{\tilde{B} b}{2(\alpha +3)} \right) |\nabla \tilde{v}_n|^2 _2 + o_n(1),
\end{align*}
which gives that $|\nabla \tilde{v}_n|_2 \to 0$. In addition, in view of $u_\rho \in S_\rho$, we have $|\tilde{v}_n|_2 \to 0$.
Thus,
\[
u_n \to u_\rho \text{ in } H_r^1(\mathbb{R}^3),
\]
and
\[
\Phi(u_\rho) = m_\rho.
\]
Finally, by Lemma 3.6 and (\ref{E310}), we deduce that $u_\rho$ solves (\ref{p})-(\ref{q}) for some $\lambda_\rho < 0$. It is easy to know that $u_\rho$ is nonnegative and radially decreasing. Moreover we also show that $u_\rho > 0$ according to the strong maximum principle.

$\hfill\Box$\\

\par\noindent{\large\bf Acknowledgements}\\
This work was supported by National Natural Science Foundation of China (Grant No.  12561023),  the Provincial Natural Science Foundation of Jiangxi (Grant No.  20252BAC250127).\\
\\

\par\noindent{\large\bf Competing interests}\\
The authors declare there is no conflict of interest exists in the submission of this manuscript, and manuscript is approved by all authors for publication.\\

\end{document}